\documentclass[11pt]{article}
\usepackage{amssymb}
\usepackage{amsmath}
\usepackage{amsthm}
\usepackage[mathscr]{eucal}
\usepackage{enumerate}

\newcommand{\V}{\mathbb{V}}
\newcommand{\U}{\mathbb{U}}
\newcommand{\R}{\mathbb{R}}

\newcommand{\N}{\mathbb{N}}

\newcommand{\Z}{\mathbb{Z}}
\newcommand{\Q}{\mathbb{Q_+}}
\newcommand{\F}{\mathbb{F}}
\newcommand{\G}{\mathbb{G}}
\newcommand{\M}{\mathbb{M}}
\newcommand{\X}{\mathbb{X}}

\theoremstyle{plain}
\newtheorem{thm}{Theorem}[section]
\newtheorem{lem}[thm]{Lemma}
\newtheorem{cor}[thm]{Corollary}

\newtheorem{proposition}[thm]{Proposition}
\theoremstyle{definition}
\newtheorem{defin}[thm]{Definition}
\newtheorem{example}[thm]{Example}
\newtheorem{remark}[thm] {Remark}

\begin{document}

\title{Divisibility of countable metric spaces} 
\author{Christian Delhomm\'e\\ E.R.M.I.T. \\ D\'epartement de Math\'ematiques et
d'Informatique\\ Universit\'e de La R\'eunion\\ 15, avenue Ren\'e
Cassin, BP 71551\\ 97715 Saint-Denis Messag. Cedex 9, La R\'eunion, France \\ {\tt
delhomme@univ-reunion.fr} 
\and Claude Laflamme\thanks{Supported by NSERC
of Canada Grant \# 690404}\\ 
University of Calgary\\
Department of Mathematics and Statistics\\ Calgary, Alberta, Canada T2N 1N4\\ 
{\tt laf@math.ucalgary.ca} 
\and Maurice Pouzet \\PCS, Universit\'e
Claude-Bernard Lyon1,\\ Domaine de Gerland -b\^at. Recherche [B], 50
avenue Tony-Garnier, \\F$69365$ Lyon cedex 07, France\\{\tt
pouzet@univ-lyon1.fr } \and Norbert Sauer\thanks{Supported by NSERC of
Canada Grant \# 691325} \\ University of Calgary \\ Department of
Mathematics and Statistics, \\ Calgary, Alberta, Canada T2N 1N4\\ {\tt
nsauer@math.ucalgary.ca} } \date{ } \maketitle

\pagebreak

\noindent \begin{tabular}{ll}
{\bf Proposed running head:} &  Divisibility of countable metric spaces \\
{\bf Contact:} & Claude  Laflamme \\
& University of Calgary\\
& Department of Mathematics and Statistics\\ 
& Calgary, Alberta, Canada T2N 1N4\\ 
& laf@math.ucalgary.ca
\end{tabular}

\pagebreak

\begin{abstract} 
Prompted by a recent question of G. Hjorth \cite{Hjorth} as to whether
a bounded Urysohn space is indivisible, that is to say has the
property that any partition into finitely many pieces has one piece
which contains an isometric copy of the space, we answer this
question and more generally investigate partitions of countable metric
spaces.

We show that an indivisible metric space must be totally Cantor
disconnected, which implies in particular that every Urysohn space $\mathbb{U}_{V}$ with $V$ 
bounded or not but dense in some initial segment of $\R_+$, is divisible.  
On the other hand we also show that one
can remove ``large'' pieces from a bounded Urysohn space with the
remainder still inducing a copy of this space, providing a certain
``measure'' of the indivisibility.  Associated with every totally
Cantor disconnected space is an ultrametric space, and we go on to
characterize the countable ultrametric spaces which are homogeneous
and indivisible.
\end{abstract}

\section{Introduction and  basic notions}

A metric space $\mathbb{M}:=(M;d)$ is called {\em divisible} if there
is a partition of $M$ into two parts, none of which contains an
isometric copy of $\mathbb{M}$. If $\mathbb{M}$ is not divisible then
it is called {\em indivisible}. Note that by repeated partition of $M$
into two pieces we obtain that if $\mathbb{M}$ is indivisible then for
every partition of $M$ into finitely many pieces there is one piece
which contains an isometric copy of the whole space. Every finite
metric space (with at least two elements ) is divisible, so the
interest lies in infinite metric spaces. The uncountable case is
different as the indivisibility property may fail badly.  For example,
every uncountable separable metric space can be divided into two parts
such that no part contains a copy of the space via a one-to-one
continuous map. This result, based on the Bernstein property 1908)(see
\cite {kuratowski} p.422) does not really involves the structure of
metric spaces.  In this paper we deal essentially with the countable
case.

After the extension of the above result to uncountable subchains of
the real line (Dushnik, Miller, 1940), the notion of indivisibility
was considered for chains and then for relational structures (see for
example \cite{Fra} \cite{pouzet}). The notion we consider also falls
under the framework of relational structures. Indeed, a metric space can
be interpreted `in several ways' to be a relational structure whose
relations are binary and symmetric, the isometries being the
isomorphisms of the relational structure. Because of this connection,
we will use some basic notions and results about relational
structures, and what we need is listed in Section~\ref{relational}.

We will show that every indivisible countable metric space is Cantor
disconnected, hence in particular, that the bounded Urysohn metric
space $\mathbb{U}_{\Q,\leq 1}$, which is Cantor connected, is
divisible.  On the other hand we will show that the space
$\mathbb{U}_{\Q,\leq 1}$ is ``almost'' indivisible, in the
sense that we can remove ``almost'' all of the elements of the space
in various ways and the remainder still contains an isometric copy of
the space.

Ultrametric spaces are special cases of totally Cantor disconnected
spaces. We will characterize the indivisible homogeneous one. It seems
to be the case that indivisible totally Cantor disconnected spaces are
rare and that there is probably no good characterization of such
spaces. We will provide various examples of indivisible countable
metric spaces.
 
\subsection {Relational structures, homogeneous structures and their ages}\label{relational}

A relational structure is a pair $\mathrm{A}:=(A; \mathbf{R})$ where
$\mathbf{R}:= (R_i)_{i\in I}$ is made of relations on the set $A$, the
relation $R_i $ being an $n_i$-ary relation identified with a subset
of $E^{n_i}$. The family $\mu:= (n_i)_{i\in I}$ is the {\it signature}
of $\mathrm{A}$. To $\mu:= (n_i)_{i\in I}$, one may attach a family
$\rho:= (r_i)_{i\in I}$ of predicate symbols and one may see
$\mathrm{A}$ as a realization of the languages whose non logical
symbols are these predicate symbols.  Let $F$ be a subset of $A$, the
induced substructure on $A$ is denoted $\mathrm A_{\restriction F
}$. Let $\mathrm{A'}:=(A'; \mathbf{R'})$ having the same signature as
$\mathrm{A}$. A local isomorphism from $\mathrm{A}$ to $\mathrm{A'}$
is an isomorphism $f$ from an induced substructure of $\mathrm{A}$
onto an induced substructure of $\mathrm{A'}$; if the domain of $f$ is
$A$ then $f$ is an {\em embedding} of $\mathrm{A}$ to
$\mathrm{A'}$. The image of an embedding of $\mathrm{A}$ in $A'$ is called a
{\em copy} of $\mathrm{A}$ in $\mathrm{A'}$.



A relational structure $\mathrm{A}:=(A; \mathbf{R})$ is {\em
divisible} if there is a partition $A=X\cup Y$ none of $X$ and $Y$
containing a copy of $\mathrm{A}$.  A relational structure which is
not divisible is called {\em indivisible}. The {\em age} of a relational
structure is the class of all finite relational structures which have
an embedding into the structure.

We will use several properties of homogeneous structures (also called
ultrahomogenous structures). Most are restatements or consequences of
the Theorem of R. Fra\"\i ss\'e (Point~\ref{frais:thm} below).
A more detailed
account can be found in the book \cite{Fra}.

\begin{enumerate}
\item A countable relational structure $\mathrm{H}:=(H,\mathbf{R})$ is
{\em homogeneous} if every local isomorphism defined on a finite
subset of ${H}$ into ${H}$ has an extension to an automorphism of
$\mathrm{H}$.

\item A countable relational structure $\mathrm{H}:=(H,\mathbf{R})$ is
homogeneous if and only if it satisfies the following {\em mapping
extension property}:
\vskip 5pt
\noindent
{\em If $\mathrm{F}:=(F;\mathbf{R})$ is an element of the age of
$\mathrm{H}$ for which the substructure of $\mathrm{H}$ induced on
$H\cap F$ is equal to the substructure of $\mathrm{F}$ induced on
$H\cap F$ then there exists an embedding of $\mathrm{F}$ into
$\mathrm{H}$ which is the identity on $H\cap F$.}

\item Two countable homogeneous structures with the same age are isomorphic.

\item A class $\mathcal D$ of relational structures has the {\it
amalgamation property} (in brief AP) if for every members $A,B,C$ of
$\mathcal D$, embeddings $f:A\rightarrow B$, $g:A\rightarrow C$, there
is some member $A'$ of $\mathcal D$ and embeddings $f': B\rightarrow
A'$, $g': C\rightarrow A'$ such that $f'\circ f=g'\circ g$.

\item\label{universality} 
A homogeneous structure embeds any countable younger structure,
i.e. any countable structure whose age is included in that of the
homogeneous one.

\item\label{frais:thm} 
A class $\mathcal D$ of finite relational structures is the age
of a countable homogeneous structure if and only if it is non-empty,
is closed under embeddability and has the amalgamation property.

\item A subset $S\not=\emptyset$ of $H$ is an {\em orbit} of
$\mathrm{H}$ if it is an orbit for the action of the automorphism
group $Aut(\mathrm{H})$ of $\mathrm{H}$ which fixes pointwise a finite
subset of $H$. That is to say that there exists a finite subset $F$ of
$H$, called a {\em socket} of the orbit $S$, so that for some $s \in H\setminus F$: 
\[ S: =\{f(s): f\in Aut(\mathrm{H}) \mbox{ and  } f(y)=y  \mbox{ for all } y \in F\}.\]

\item If $\mathrm{H}$ is a countable homogeneous structure, then a
subset $S \subseteq H$ is an orbit of $\mathrm{H}$ if there is an $s
\in H \setminus F$ and $S$ is equal to the set of all elements $t\in
H$ so that the function which fixes the socket $F$ pointwise and maps
$s$ to $t$ is an isomorphism of the substructure of $\mathrm{H}$
induced on $S\cup\{s\}$ on the substructure of $\mathrm{H}$ induced on
$S\cup\{t\}$.  That is, the orbit $S$ is the set of all elements of
$H$ which are of the same ``one-type'' over $F$.

\item\label{prolong:orbit}
If $\mathrm{H}$ is a countable homogeneous structure, a subset
$X\subseteq H$ induces an isomorphic copy of $\mathrm{H}$ if and only
if $S\cap X\not= \emptyset$ for every orbit $S$ of $\mathrm{H}$ with 
socket  a subset of $X$.

\item\label{kappa:coloring}
 Let $\kappa$ be a cardinal and
$\mathcal{A_{\kappa}}$(resp. $\mathcal{A_{ \kappa, <\omega,}}$) be the
collection of all (resp. finite) relational structures
$\mathrm{B}:=(B; \mathbf{R})$ where $\mathbf{R}:=(R_i)_{i<\kappa}$ is
a sequence of irreflexive and symmetric binary relations symbols for which for
all $x,y\in B$ with $x\not= y$ there exists exactly one $i<\kappa$
with $R_i(x,y)$.  The class $\mathcal{A_{ \kappa, <\omega}}$ has the
amalgamation property.  If $\kappa\leq \omega$ then it is countable,
therefore, this is the age of a countable homogeneous structure, that
we denote $\mathrm{H}_{\kappa}$.  For example, $\mathrm{H}_2$ is the
well-known {\it Random graph} or { \it Rado graph}. Each such $H_{k}$  is indivisible.
\end{enumerate}

\subsection{Metric spaces and relational structures}

Let us recall a few standard notions. Given two metric spaces
$\mathbb{M}:=(M;d)$ and $\mathbb{M'}:=(M';d')$, a {\it local isometry}
from $\mathbb{M}$ to $\mathbb{M'}$ is an isometry $f$ from a subspace
of $\mathbb{M}$ onto a subspace of $\mathbb{M}$, and this is an {\it
isometric embedding} if the domain of $f$ is $M$. \\ $\mathbb{M}$ is
called {\it homogeneous } if every local isometry defined on $M$ and
with values in $M$ extends to an isometry from $\mathbb{M}$ onto
itself. \\ The {\it age} of $\mathbb{M}$ is the collection of finite
metric spaces which embed into $\mathbb{M}$. \\ Finally the {\em
spectrum } of $a\in M$ is the set $Spec(\mathbb{M}, a)=\{d(a,x) \mid
x\in M\}$ and the {\em spectrum} of $(M;d)$ is the set
$Spec(\mathbb{M}):=\cup \{Spec(\mathbb{M}, a): a\in M\}= \{d(x,y) \mid
x,y\in M\}$.

Metric spaces also fall under the realm of relational structures in
various ways.  To exemplify this association, consider a set $I$, a
map $f:I\rightarrow \R_+$ and set $\mu:= (n_i)_{i\in I}$, where
$n_i:=2$ for all $i\in I$.  To a metric space $\mathbb{M}:=(M;d)$
associate two relational structures, namely $\mathbb{M}_{f, \leq}:=(M;
\mathbf {R})$ and $\mathbb{M}_{f,=}:=(M; \mathbf {S})$ where $\mathbf
{R}:= (R_i)_{i\in I}$ and $\mathbf {S}:= (S_i)_{i\in I}$ are defined
by:
\begin{equation}\label{metric1} 
(x,y)\in R_{i}\Longleftrightarrow d(x,y)\leq f(i)
\end{equation}
\begin{equation}\label{metric2}
(x,y)\in S_{i}\Longleftrightarrow 0\not =d(x,y)=f(i)
\end{equation}

Using the above notation, the following result summarizes the
connections we will need, and the straightforward proof is left to the
reader.

\begin{lem} \begin{enumerate}
\item Every local isometry of $\mathbb{M}$ is a local isomorphism of
$\mathbb{M}_{f, \leq}$ and of $\mathbb{M}_{f, =}$.

\item Every local isomorphism of $\mathbb{M}_{f, \leq}$ (resp. of
$\mathbb{M}_{f,=}$) is a local isometry of $\mathbb{M}$ if and only if
:
\begin{enumerate} 
\item either the spectrum of $\mathbb{M}$ contains at most a non-zero
element,
\item or the image of $f$ separates the spectrum of $\mathbb{M}$ in
the sense that $f(I)\cap [p,p')\not=\emptyset $, (resp in the sense
that $f(I)\cap \{p,p'\}\not= \emptyset$) for every $p,p'\in
Spec(\mathbb {M})$ such that $0<p<p'$.

\end{enumerate}
\end{enumerate} 
\end{lem}

Conversely, every binary relational structure $\mathrm{B}:=(B;
\mathbf{R})$ can be viewed as a metric space, provided that the number
of isomorphic types of induced substructures on two element subsets of
$B$ is not greater than the continuum.  Indeed, let $a\in
\R_+\setminus \{0\}$ be given, we may define a one-to-one map
$\varphi: [B]^2\rightarrow [a, 2a]$ such that
$\varphi(\{x,y\})=\varphi(\{x',y'\})$ if and only if $\mathrm
B_{\restriction \{x,y\}}$ and $\mathrm B_{\restriction \{x',y'\}}$ are
isomorphic. The map $d: B\times B \rightarrow [a, 2a]$ defined by
setting $d(x,y):=\varphi(\{x,y \})$ if $x\not =y$ and $d(x,y):=0$
if $x =y$ is a distance. This is particularly the case if
$\mathrm{B}\in \mathcal A_{\kappa}$. A map $f: \kappa\rightarrow [a,
2a]$ be given, we will set $m_f(\mathrm {B}):=(B, d)$ where
$d(x,y):=f(i)$ if $R_i(x,y)$ and $d(x,y):=0$ otherwise (if $x=y$).

\subsection{Homogeneous metric spaces}

Let $V$ be a set such that $0\in V\subseteq \R_+$. Let $\mathcal
M_V$ (resp. $\mathcal M_{V, <\omega}$) be the collection of metric
spaces (resp. finite metric spaces) $\mathbb{M}$ whose spectrum is
included into $V$. We may note that any such $V$ is in fact a
spectrum,  indeed $V= Spec (\mathbb{M})$ where $\mathbb{M}:=(V;d)$ and
$d(x,y):= \max (\{x,y\})$.

For  $u_1$, $u_2$, $u_1'$ , $u_2'$ in $V$, let  
\[
\phi(u_1,u_2,u_1',u_2'):=[\max \{|u_1-u_2|,|u_1'-u_2'|\},\min \{u_1+u_2,u_1'+u_2'\}]
\]
and set 
\[
\rho_V(u_1,u_2,u_1',u_2')\;  \text{ if} \; \phi(u_1,u_2,u_1',u_2')\cap V\neq\varnothing
\]
We say that $V$ satisfies the {\it four-values condition} if 
\[
\forall u_1,u_2,u_1',u_2'\text{ in }V\
\bigl(
\rho_V(u_1,u_2,u_1',u_2')\Rightarrow\rho_V(u_1,u_1',u_2,u_2')
\bigr)
\]

\begin{lem}\label{lem:4vc-hered}
Let   $u_1$, $u_2$, $u_1', u_2'\in V$
such that  $\rho_V(u_1,u_2,u_1',u_2')$ holds. 
\begin{enumerate}
  \item If  $u_1$, $u_2$, $u_1'$ and  $u_2'$ are all non zero then 
   $\phi(u_1,u_2,u_1',u_2')\cap(V\setminus\{0\})\neq\varnothing$.
  \item If some argument  is zero,
then $\rho_V(u_1,u_1',u_2,u_2')$ holds.
In particular the four-value condition is equivalent to its restriction to non zero arguments.
\end{enumerate} 
\end{lem}

\begin{proof}
\begin{enumerate}
  \item If $\max \{|u_1-u_2|,|u_1'-u_2'|\}>0$ then every element of $\phi(u_1,u_2,u_1',u_2')$
is positive.
If $\max \{|u_1-u_2|,|u_1'-u_2'|\}=0$ then $\min \{u_1,u_2,u_1',u_2'\}\in\phi(u_1,u_2,u_1',u_2')\cap(V\setminus\{0\})$.
  \item We may assume wlog that $u_1=0$. Since $\rho_V(u_1,u_2,u_1',u_2')$ holds,
  $\phi(u_1,u_2,u_1',u_2')=\{u_2\}$
and likewise $\phi(u_1,u_1',u_2,u_2')=\{u_1'\}$.
\end{enumerate}
\end{proof}

\begin{proposition}Let $V$ be a set such that $0\in V\subseteq \R_+$.
The class $\mathcal M_{V, <\omega}$ is the age of a metric space whose
spectrum is $V$.
Furthermore the following are equivalent~:
\begin{enumerate}
  \item  $\mathcal M_{V, <\omega}$ has the amalgamation property;
   \item
   $\mathcal M_{V, <\omega}$ has the {\em disjoint amalgamation property}~:
   i.e. two members of $\mathcal M_{V, <\omega}$ that coincide on their intersection admit (on their union) a common extension in $\mathcal M_{V, <\omega}$;
\item 
For any two members $(M_1,d_1)$ and $(M_2,d_2)$ of $\mathcal M_{V, <\omega}$
such that $d_1$ and $d_2$ coincide on $M_1\cap M_2$ and 
such that $|M_1|=|M_2|=3$ and $|M_1\cap M_2|=2$,
there is a semi-distance on $M_2\cup M_2$ with spectrum included in $V$ and whose restrictions to $M_1$ and $M_2$ are $d_1$ and $d_2$ respectively.
(A semi-distance is symmetric and satisfies the triangular inequality but may fail to satisfy the separation condition.)
\item $V$ satisfies the four-values condition.
\end{enumerate}
\end{proposition}

\begin{proof}[Proof of the proposition] 
First, there is a family $(\mathbb {M}_i)_{i\in I}$ of at most
$\kappa:=\vert V\vert \cdot {\aleph_{0}}$ members of $\mathcal M_{V,
<\omega}$ such that every member of $\mathcal M_{V, <\omega}$ embeds
into one of the $\mathbb {M}_i$'s. Pick an element $0_i\in M_i$ for
each $i\in I$. Set $M:= \{\overline x:= (x_i)_{i\in I}: x_i \in M_i$
for all $i\in I$ and $\{i\in I: x_i\not = 0_i\}$ is finite $\}$, set
$d(\overline x, \overline y):= \max \{d(x_i,y_i): i\in I\}$ and set
$\mathbb M:= (M, d)$. Then the age of $\mathbb M$ is $\mathcal M_{V,
<\omega}$. Since every subset of $\R_+$ containing $0$ is a spectrum,
the spectrum of $\mathbb M$ is $V$.

Next, the implications  $2\Rightarrow 1\Rightarrow 3$ are obvious.  We prove  $3\Rightarrow 4\Rightarrow 2$.

\noindent$\bullet$
We assume that Point~3 holds and we check the four-value condition~:
Consider $u_1$, $u_2$, $u_1'$ and $u_2'$ in $V$
and assume that $\rho_V(u_1,u_2,u_1',u_2')$ holds.
Thanks to the second part of the Lemma, 
we can assume that $u_1$, $u_2$, $u_1'$ and $u_2'$ 
are all positive, and then thanks to the first part, 
we know that $\phi(u_1,u_2,u_1',u_2')$
contains some non zero element $v$ of $V$.
Then given a four-element set $\{x_1,x_2,y,y'\}$,
the following define distances $d_1$ and $d_2$ on $M_1:=\{x_1,y,y'\}$
and on $M_2:=\{x_2,y,y'\}$ that coincide on $M_1\cap M_2=\{y,y'\}$~:
 \[
\begin{cases}
d_1(x_1,y)=u_1,\ d_1(x_1,y')=u_1',\ d_1(y,y')=v\\
d_2(x_2,y)=u_2,\ d_2(x_2,y')=u_2',\ d_2(y,y')=v
\end{cases}
\]
Now by Point~3, $d_1$ and $d_2$ have a common extension to a semi-distance $d$ on $M_1\cup M_2$ with values in $V$.
Then $d(x_1,x_2)$ belong to $\phi(u_1,u_2,u_1',u_2')\cap V$.

\noindent$\bullet$
Now we assume that the four-value condition holds and we show the disjoint amalgamation property~:
Consider two members $\mathbb M_1=(M_1,d_1)$ and $\mathbb M_2=(M_2,d_2)$ of $\mathcal M_{V, <\omega}$
such that $d_1$ and $d_2$ coincide on $M_1\cap M_2$.,

 First consider the case where both $M_1\setminus M_2$
 and $M_2\setminus M_1$ are singletons; then let $x_1$ and $x_2$ denote their
 respective elements.
Observe that extending $d_1$ and $d_2$ by setting $d(x_1,x_2)=d(x_2,x_2)=w$ yields a distance extending $d_1$ and $d_2$
if and only $w>0$ and the triangular inequalities involving both $x_1$ and $x_2$ are 
satisfied (the other ones involve only one of $d_1$ or $d_2$), and that holds 
precisely when 
\[
\forall z\in M_1\cap M_2\ 
\vert d_1(x_1,z)-d_2(x_2,z) \vert
\leq w\leq
d_1(x_1,z)+d_2(x_2, z)
\]
 Besides, it follows from the triangular inequality that
 $a':=\max \{\vert d_1(x_1,z)-d_2(x_2,z) \vert : z\in M_1\cap M_2\}\leq a:=\min  \{d_1(x_1,z)+d_2(x_2, z): z\in M_1\cap M_2\}$.
 Then pick $y$ and $y'$ in $M_1\cap M_2$ such that 
 $|d_1(x_1,y')-d_2(x_2,y')|=a'$
 and $d_1(x_1,y)+d_2(x_2, y)=a$, and let
 \[
\begin{cases}
u_1:=d_1(x_1,y),\ u_1':=d_1(x_1,y')\\
u_2:=d_2(x_2,y),\ u_2':=d_2(x_2,y')
\end{cases}
\]
All those values are positive members of $V$
and the distance between $y$ and $y'$ attests that 
$\rho_V(u_1,u_2,u_1',u_2')$ holds,
thus by the four-value condition and the first part of the Lemma,
and given that $\phi(u_1,u_1',u_2,u_2')=[a,b]$,
the intersection $[a',a]\cap V$ is non-empty.
That concludes that case.

 Now we proceed by induction on the cardinality $m$ of the symmetric difference
 $M_1\vartriangle M_2$.
 If $m\leq 1$ then $M_1\subseteq M_2$ or $M_2\subseteq M_1$,
 in which case there is nothing to prove.
 So assume that $m>1$ and neither $M_1\subseteq M_2$ nor $M_2\subseteq M_1$.
 Pick any $x_1\in M_1\setminus M_2$ and $x_2\in M_2\setminus M_1$.
 Let $M:=M_1\cup M_2$.
 By induction assumption,
 $\mathbb{M}_1$ and $\mathbb{M}_2\restriction(M_2\setminus\{x_2\})$
 admit a common extension $\mathbb{M}_1'$ on $M_1\cup(M_2\setminus\{x_2\})=M\setminus\{x_2\}$
 and then, still by induction assumption
 $\mathbb{M}_1'$ and $\mathbb{M}_2$
 admit a common extension $\mathbb{M}_2'$ on $(M\setminus\{x_2\})\cup M_2=M\setminus\{x_1\}$.
 Then by the first case, $\mathbb{M}_1'$ and $\mathbb{M}_2'$,
 which coincide on $M\setminus\{x_1,x_2\}$,
admit a common extension on $M$,
 which extension then extends $\mathbb{M}_1$ and $\mathbb{M}_2$.
\end{proof}

Fra\"\i ss\'e's theorem  (Point~\ref{frais:thm} above) gives immediately:

\begin{thm}
If $V$ is countable then it  satisfies the four-values condition if and only if  there is a countable homogeneous space
$\mathbb{U}_V$ whose age is $\mathcal M_{V, <\omega}$.
\end{thm}

We call the space $\mathbb{U}_V$ the {\em Urysohn space with spectrum
$V$}. If $V:=\Q$ then $\mathbb{U}_{\mathbb Q_+}$ is the homogeneous metric space whose age
is the set of all finite metric spaces whose spectrum is a subset of
the set of rationals. The Cauchy completion of $\mathbb{U}_{\mathbb{Q_+}}$ is the
famous space discovered by Urysohn \cite{urysohn}.

\begin{example}
\begin{enumerate}
   \item Suppose that for some $a\in \R_+\setminus \{0\}$, $V\setminus
\{0\}\subseteq [a, 2a]$. Then $V$ satisfies the four-values condition and in fact $\U_V= m_f (\mathrm {H}_{\kappa})$ where
$\kappa:= \vert V\setminus\{0\}\vert $ and $f: \kappa \rightarrow
V\setminus\{0\}$ is a bijective map (see Point~\ref{kappa:coloring} of Section~\ref{relational}).
 \item An example of finite set $V$ which does not satisfies the four point conditioon is 
  $V=\{0,1,3,4,5\}$.
Indeed,  $M_1:=\{y,y',x_1\}$ and $M_2:=\{y,y',x_2\}$
 with $v:=d(y,y')=4$, 
 $u_1:=d(x_1,y)=1$, $u_2:=d(x_2,y)=1$, 
  $u_1':=d(x_1,y')=5$ and $u_2':=d(x_2,y')=3$.
Each "distance" in $M_1$ (namely 1,4 and 3)  or in $M_2$ (1,4,5)
 is less than or equal to the sum of the other two, 
 thus $d$ is indeed a distance on $M_1$ and on $M_2$.
The only possible value for $d(x_1,x_2)$ 
in order that $d$ be a semi-distance on $M_1\cup M_2$ is 2~:
Indeed
$2=5-3=d(x_1,y')-d(x_2,y)\leq d(x_1,x_2)\leq d(x_1,y)+d(x_2,y)=1+1=2$.
If $V$ is infinite  is not sufficient that $V$ be dense in $]\inf V,\sup V[$ to have the
four-values conditions, since the example above also shows that $\R_+\setminus\{2\}$
fails to have it.
\item
A sufficient condition for the four-value condition is 
\[
\forall\ u_1,u_2,u_1',u_2'\in V\
u_1-u_1'\leq u_2+u_2
'\Rightarrow \exists v\in V\ 
u_1-u_1'\leq v\leq u_2+u_2
\] 
For  example,  the set $V$ of positive powers of $\frac 12$ satisfies it. 
Indeed, if
$0<\frac 1{2^p}-\frac 1{2^q}\leq\frac 1{2^r}+\frac 1{2^s}$, then $\frac 1{2^p}$
lies in between.
 
This condition is not necessary~: consider $V:=\{0,1,3,5\}$. 
 
Notice that this sufficient condition holds  
whenever $V$ is closed
under sum or  absolute value of the difference,
or more generally
when for all $a$ and $b$ in $V$, if $a+b<\sup V$ then $a+b$ belongs to $V$
  
Examples, like $\N$, $\Q$, $\{0,\dots n\}$ and their  corresponding  Urysohn spaces are considered in~\cite {pouzroux}. 
\end{enumerate}
\end{example}

\begin{lem} 
If $V$ satisfies the four-values condition then every initial segment $V'$ of $V$ satisfies this condition
\end{lem}

\begin{proof} Let $u_1,u_2, u_3,u_4 \in V'$ such that $\rho_{V'}(u_1,u_2, u_3,u_4)$ holds. Let $w\in \phi(u_1,u_2, u_3,u_4)\cap V'$ and let $r:= \max \{u_1,u_2,u_3,u_4, w\}$.
Case 1: $\min  \{u_2+u_3, u_4+u_1\}\leq r$. In this case, since $r\in V'$, $\min  \{u_2+u_3, u_4+u_1\}\in V'$. Since   $\rho_V(u_2, u_3, u_4,u_1)$ holds,  $\rho_{V'}(u_2, u_3, u_4,u_1)$  holds too.
Case 2: $r<\min  \{u_2+u_3, u_4+u_1\}$. In this case $r\in \phi(u_2,u_3, u_4,u_1)\cap V'$ thus $\rho_{V'}(u_2, u_3, u_4,u_1)$  holds.

\end{proof}

\begin{remark}
For a more intuitive  proof,  based  on the amalgamation property  of $\mathcal M_{V,<\omega}$,
let $\mathbb M_1, \mathbb M_2\in \mathcal M_{V',<\omega}$ and 
 let $\mathbb M:=(M,d)\in \mathcal M_{V,<\omega}$ be a common extension. Set $d':=d\wedge \delta$, where 
$\delta$ is the maximum of the diameters of  $\mathbb M_1$ and $\mathbb M_2$. \end{remark}

With this lemma and the above theorem follows that if there is an Urysohn space with spectrum $V$ then for every $\ell \in \R_+$ there is an Urysohn space with spectrum $V\cap [0,\ell]$. We denote this space  $\mathbb
{U}_{V,\leq \ell}$, eg $\mathbb{U}_{\Q, \leq \ell}$ is the the homogeneous metric space whose age is the set of all finite metric spaces whose spectrum is a subset of
the set of rationals  in the interval
$[0,\ell]$).

Let $V$ be  subset of $\R_+$ containing $0$; we say that $V$ is {\it residuated}   if for every  $x, y\in V$ with $x\leq y$ the set  $\{ r \in V: y\leq x+ r\}$ has a least element, denoted $y\setminus x$.  This is the case if $V$ is finite, if $V$ is the positive part of a additive subgroup of $\R$  or if  $V$ is {\it meet-closed } in the sense that for every non-empty subset of $V$ its infimum in $\R$ belongs to $V$.

The following proposition shows that the four-values condition is just what is needed in order to extend to metric spaces over $V$ the most fundamental property of ordinary metric spaces.

\begin{proposition}\label{fundamental} Let $V$ be  subset of $\R_+$ containing $0$. If $V$ is residuated, then $V$ satisfies the  four-values condition if and only if  the map: $d_V: V\times V\rightarrow V$ defined by $d_V(x,y):= \max  \{ y\setminus x, x\setminus y\}$ is a distance on $V$. When this condition is realized,  $d(x,y)= Sup\{d_V(d(x,z), d(y,z)): z\in M\}$ for every $ x, y \in \M:= (M,d)\in \mathcal M_V$. In particular, if $\M$ is finite  the map $\overline d: M\rightarrow V^M$ defined by setting $\overline d(x)(y):= d(x,y)$ is an isometric embedding of $\M$ into $V^M$ equipped with the "Sup" distance.

\end{proposition}

\begin{proof} Clearly, $d_V(x,y)=Inf\{r\in V: \vert x-y\vert \leq r\}$. It follows that $d_V$ is symmetric, $d_V(x,y)$=0 iff $x=y$ and $d_V(0, x)=x$ for every $x\in V$. Suppose that the four-values condition holds. Let $x,y, z\in V$; set $u_1:= x,u_2:= d_V(x,z), u_3:=d_V(z,y), u_4:= y$. We have $z\in \phi(u_1,u_2, u_3,u_4)\cap V$, proving that  $\rho_V(u_1,u_2,u_3,u_4)$ holds. Since the four-values condition holds, $\rho_V(u_2,u_3,u_4,u_5) $ holds, that is  $\max  \{\vert u_2-u_3\vert, \vert u_4-u_1\vert \} \leq r \leq  \min  \{u_2+u_3, u_4+u_1\}$ for some $r\in V$. In particular $\vert u_4-u_1\vert \leq r \leq u_2+u_3$ that is $\vert x-y\vert \leq r \leq  d_V(x,z)+d_V(z,y)$. The triangular inequality  $d_V(x,y)\leq d_V(x,z)+d_V(z,y)$ follows. Conversely, suppose that $d_V$ is a distance on $V$. Let $u_1,u_2, u_3,u_4 \in V$ such that $\rho_V(u_1,u_2, u_3,u_4)$ holds. Let $r:= \max \{d_V(u_2,u_3), d_V(u_4,u_1) \}$.  First, we have  $r\in V$; next, from the definition of $d_V$,  we have $\vert u_2-u_3\vert \leq d_V(u_2,u_3)$ and $\vert u_4-u_1\vert \leq d_V(u_4,u_1)$ hence $\max  \{\vert u_2-u_3\vert, \vert u_4-u_1\vert\} \leq r$;  finally   the triangular inequality   $d_V(u_4, u_1)\leq d_V(u_4, w)+d(w,u_1)$ applied first  to $w:=0$ gives $d_V(u_4,u_1)\leq  u_4+u_1$ and applied to  $w\in \phi(u_1,u_2,u_3,u_4)\cap V$ gives  $d_V(u_4, u_1)\leq u_3+u_2$ hence $d_V(u_4,u_1)\leq \min \{u_4+u_1,  u_3+u_2\}$; the same gives $d_V(u_2,u_3)\leq \min \{u_4+u_1,  u_3+u_2\}$, hence $r\leq \min \{u_4+u_1,  u_3+u_2\}$. This proves that  $r\in \phi(u_2, u_3, u_4,u_1)\cap V$, hence $\rho _V(u_2, u_3, u_4,u_1)$ holds. Thus the four-values condition holds.
If $\M:=(M,d)$ is an arbitrary metric space with spectrum included in $V$, the equality $d(x,y)= Sup\{d_V(d(x,z), d(y,z)): z\in M\}$ is obvious. If $M$ is finite, it  expresses  the fact that $\overline d$ is an isometric embedding from $\M$ into $V^M$ equipped with the Sup-distance. 
\end{proof}

\begin{remark} As it is well-known every  metric space embeds into some $\ell^{\infty}_{\R}$ space equipped with the Sup distance. A similar result holds for members of $\mathcal M_V$ provided that    $V$  is meet closed and satisfies the four-values condition. 
\end{remark}


As we shall see later on, some $\U_V$'s are divisible, still
all metric spaces with age $\mathcal M_{V, <\omega}$ satisfy a weaker
version of the indivisibility property~:

\begin{thm} \label {general}
Let $\mathbb {M}$ be a metric space with age $\mathcal M_{V,
<\omega}$.

\begin{enumerate}

\item For every partition of $M$ into two parts $X$ and $Y$ one of the
induced metric spaces $\mathbb M_{\restriction X}$ and $\mathbb
M_{\restriction Y}$ has the same age as $\mathcal M_{V, <\omega}$.

\item If $V$ is countable and bounded then there is an indivisible
metric space with age $\mathcal M_{V, <\omega}$

\end{enumerate}
\end{thm}

\begin{proof}
This result is due to the fact that $\mathcal M_{V, <\omega}$ is
closed under finite product. This is a special case of Corollary 1 of \cite {pouzros} built on  \cite {jezek}.  The key ingredient is the
Hales-Jewett 's theorem \cite {halesjewett}.
 
{\bf Claim} For every $\mathbb{F}\in \mathcal M_{V, <\omega}$ there is
some $\mathbb{G}\in \mathcal M_{V, <\omega}$ such that for every
partition of $G$ into two parts $X$ and $Y$, one of the spaces induced by
$\mathbb{G}$ embeds $\mathbb{F}$.  Recall that a {\it combinatorial
line} of a finite cartesian power $N^{n}$ of $N$ is a set of the form
$L_{\overline l}:= \{\overline x:= (x_i)_{i<n} \in N^n: x_i= l_i $ for
all $i\in K$ and $x_i=x_j$ for all $i,j\not \in K\}$ where $l:=
(l_i)_{i\in K}$ and $K\subset n$.  According to Hales-Jewett 's
theorem, if $n$ is large enough then for every partition of $N^{n}$
into two parts $X$, $Y$ one of the parts contains a combinatorial
line. Thus, if we equip $F^n$ with the "sup-distance", the resulting
space $\mathbb{G'}$ satisfies the conclusion of the claim.

To prove part $1$, let $\mathbb {M}$ be metric space with age
$\mathcal M_{V, <\omega}$ and let $X, Y$ be a partition of $M$. Assume
for a contradiction that the ages $\mathcal A$ of ${\mathbb
M}_{\restriction X}$ and $\mathcal B$ of ${\mathbb M}_{\restriction
Y}$ are distinct from $\mathcal M_{V, <\omega}$, and thus let $\mathbb
{M}_X\in \mathcal M_{V, <\omega}\setminus \mathcal A$ and $\mathbb
{M}_Y\in \mathcal M_{V, <\omega}\setminus \mathcal B$.  Select $A, B
\subseteq M$ such that $\mathbb M_{\restriction A}$ and $\mathbb
M_{\restriction B}$ are an isometric copy of $\mathbb {M}_X$ and
$\mathbb {M}_Y$ respectively.  For $\mathbb
F:=\mathbb{M}_{\restriction A\cup B}$ there is no $\mathbb {G}$
satisfying the conclusion of the claim, a contradiction.

Now to prove $2$, let $a\in V$ such that $2a$ is an upper-bound of
$V$.  Let $(\F_n)_ {n<\omega}$ be an enumeration of the members of
$\mathcal M_{V, <\omega}$. According to the Claim above, there is a
sequence $(\G_n)_{n<\omega}$ such that $\G_{n+1}$ contains an
isometric copy of $\F_{n+1}$ and for every partition of $G_{n+1}$ into
two parts one of the part contains an isometric copy of $\G_n$. Let
$G$ be the disjoint union of the $G_n$'s and $d: G\times G\rightarrow
V$ be defined by $d(x,y):= d_n(x,y)$ if $x,y\in G_n$ for some $n$ and
$d(x,y):= a$. Then $\G:= (G,d)$ is an indivisible metric space with
age $\mathcal M_{V, <\omega}$.
\end{proof}

Theorem~\ref{thm:unbounded} asserts that the condition that $V$ is
bounded is necessary.

\subsection{Indivisibility of Urysohn spaces}

Here is a short summary of indivisibility results regarding the
Urysohn spaces.

\noindent
$\bullet $ $\mathbb U_V$ is indivisible if $V\subseteq \{0\}\cup [a,
2a]$ for some $a>0$. Indeed in this case $\mathbb U_V= m_f(\mathrm
{H}_{\kappa})$, $\kappa\leq \omega$.  Since $\mathrm H_{\kappa}$ is
indivisible, it follows that $\mathbb U_V$ is indivisible as well.


\noindent $\bullet$ Let $R$ and $S$ be two relational structures with
the same signature.  Write $R\preceq S$ if there exists a partition of
$R$ into finitely many parts $R_0, R_1,\dots, R_{n-1}$ so that for all
$i\in n$ there is an embedding of $R_i$ into $S$. A necessary
condition for a homogeneous structure to be indivisible is that any
two orbits of it are related under $\preceq$, see \cite{ES}. Urysohn
metric spaces satisfy this necessary condition according to Corollary
\ref{norberttest}.  This then implies together with Item~1 of Theorem
\ref{general}, that if a homogeneous metric space is indivisible then
the ages of any two orbits are comparable under $\subseteq$.

It follows from results in \cite{ES} and \cite{Sauer} that homogeneous
binary structures with finite signature whose age has free
amalgamation are indivisible if and only if they satisfy that
necessary condition above. (It seems to one of us, Sauer, that this
result could be extended without too much of a problem to homogeneous
structures with free amalgamation and infinite binary signature.)

Ages of homogeneous metric spaces satisfy the weaker notion of strong
amalgamation. (See the appendix of \cite{Fra} for the definitions of
free and strong amalgamation.) The Urysohn space
$\mathbb{U}_{\Q,\leq 1}$, and of course by a similar argument
the Urysohn space $\mathbb{U}_{\Q,\leq a}$ for every positive
real $a$, is divisible according to Theorem~\ref{thm:Urdis}; providing
a large number of examples of homogeneous structures which satisfy the
neccessary condition above and which have strong amalgamation but
which are divisible.

If $V:=\{0,1\}$ the metric space $\mathbb U_V$ is indivisible. If
$V:=\{0,1,2\}$ then the homogeneous metric space $\mathbb U_V$ is just
a cryptomorphic version of the Rado graph.  Associate with every edge
of the Rado graph distance 1 and with every non-edge distance 2. Hence
it follows that $\mathbb U_V$ is indivisible in this case.

Let $V:= \{0,1,2,3\}$. Then $\mathbb U_V$ is a cryptomorphic version
of the homogeneous graph $H$ with two types of edges, $E_1$ and $E_3$,
which does not contain a triangle with two edges of type $E_1$ and one
edge of type $E_3$. Associate with every edge of type $E_1$ distance 1
and with every edge of type $E_3$ distance 3 and with every non-edge
distance 2. The homogeneous structure $H$ has free amalgamation and
satisfies the chain condition. Hence it follows from Corollary 8.2 of
\cite{Sauer} that $H$ and therefore $\mathbb U_V$ is indivisible.

We do not know if $\mathbb U_V$ is indivisible for $V:=\{0,1,2,3,4\}$.

More generally the situation is as follows. Let $V$ with $0\in V$ be a
finite set of non negative real numbers satisfying the conditions that
there is a number $0\not= a\in V$ so that:
\begin{align}\label{three}
2\cdot\min(V) \leq a \text{ and } \min(V)+a\leq \max(V).
\end{align}
It follows that in this case the conditions of Lemma 1.2 are satisfied. Hence the age $\mathcal{M}_{V,<\omega}$ has amalgamation and there exists a homogeneous metric space $\mathbb{U}_V$.

This metric space $\mathbb{U}_V$ is a cryptomorphic version of the
homogeneous graph $H$ with several types of edges $E_i$ for $i\in
V\setminus\{0,a\}$. If we associate in $H$ with the edges $E_i$ the
distance $i$ and with the non-edge between different vertices the
distance $a$ we will obtain the metric space $\mathbb{U}_V$. The graph
$H$ satisfies the chain condition because the metric space
$\mathbb{U}_V$ satisfies the chain condition. Hence, according to
Corollary 8.2 of \cite{Sauer}, $H$ is indivisible and so is
$\mathbb{U}_V$.

It follows from the definition of free amalgamation of relational
structures that the age of the graph $H$ obtained as above from the
metric space $\mathbb{U}_V$ has free amalgamtion if and only if $V$
satisfies the Inequalities~\ref{three}.

\noindent $\bullet$ If $\mathbb U_V$ is indivisible then for every
$a\in \R_+\setminus \{0\}$, $[0, a]\cap V$ is not dense in $[0,
a]$(Theorem~\ref {cantorconnected}). In particular neither $\mathbb
U_{\Q}$ nor $\mathbb{U}_{\Q,\leq \ell}$ is indivisible.

\noindent $\bullet$ If $\mathbb U_V$ is indivisible then $V$ must be
bounded (Theorem\ref {thm:unbounded}). In this case let $a:= Sup V$, must then $V \cap [0,
\frac{a}{2}]$ be dually well-founded?

\noindent $\bullet$ If $\mathbb U_V$ is indivisible then there is a map $u:
V\rightarrow \R_+$ such that $u(v)\leq v$ for every $v\in V$ and
$d^\ast:= u\circ d$ is an ultrametric distance on $U$. It follows that
$\U_V^\ast:= (U, d^\ast)$ is homogeneous and indivisible (see Theorem
\ref{thm:ultra}).




\section{Ultrametric  spaces and homogeneous ultrametric spaces}

A metric space is an {\em ultrametric space} if it satisfies the
strong triangle inequality $d(x,z)\leq \max\{d(x,y),d(y,z)\}$. See
\cite{Lemin} for example. Note that a space is an ultrametric space if
and only if $d(x,y)\geq d(y,z)\geq d(x,z)$ implies $d(x,y)=d(y,z)$.
What we did in the previous section for general metric spaces work for
ultrametric spaces.

Let $V$ be a set such that $0\in V\subseteq \R_+$. Let $\mathcal
Mult _V$ (resp. $\mathcal Mult_{V, <\omega}$) be the collection of
ultrametric metric spaces (resp. finite ultrametric spaces)
$\mathbb{M}$ whose spectrum is included into $V$. Then $\mathcal Mult
_{V, <\omega}$ is the age of a metric space whose spectrum is $V$; it
is closed under embeddability and has the amalgamation property. If
$V$ is countable then there is a countable homogeneous ultrametric
space $\mathbb{U}ult _V$ whose age is $\mathcal Mult _{V, <\omega}$
and has spectrum $V$; we call it the {\em Urysohn ultrametric space
with spectrum $V$}. We give a description of this space in Proposition \ref {homogeneous2}. 

For a given set $V$, $\U_V$ and $\mathbb{U}ult_V$ are in general
different, except if $V=\{a_n: n\in D\}$ where $D$ is an interval of
the set $\Z$ of integers and $2a_{i+1}<a_i$ for all $i, i+1 \in D$.

Homogeneous ultrametric are easy to describe. In fact ultrametric
spaces can be described by means of real-valued trees.  An ordered set
$P$ is a {\it forest} if for every $x\in P$ the set $\downarrow x:=
\{y\in P: y\leq x\}$ is a chain; this is a {\it tree} if in addition
every pair of elements of $P$ has a lower bound.  If every pair $x,
y\in P$ has an infimum, denoted $x\wedge y$, we will say that $P$ is a
{\it meet-tree}. We say that $P$ is {\it ramified} if for every $x,
y\in P$ such that $x<y$ there is some $y'\in P$ such that $x<y'$ and
$y'$ incomparable to $y$. In the sequel, we consider ramified
meet-trees such that every element is below some maximal
element. These posets are meet-semilattices generated by their
coatoms.  We will need the following property

\begin{lem} 
Let $P$ be a ramified meet-tree such that every element is below some
maximal element.  For every $x\in P\setminus \max (P)$ there is a subset
$X\subseteq \max  (P)$ of maximum cardinality such that $x=a\wedge b$
for every pair of distinct elements $a, b$ of $X$
\end{lem}

\begin{proof}
For two elements $a$ and $b$ above an element $x$, set $a\equiv b$ if
$x\not = a\wedge b$. Observe that this is an equivalence relation. A
set $X$ which meets each equivalence classe has maximum size.
\end{proof}

The cardinality of $X$, denoted $d_P(x)$, is the {\it degree} of $x$.
For $x\in \max (P)$ we set $d_P(x):= 0$. If $P$ is finite or
well-founded, this the ordinary notion of degree, that is the number
of upper-covers of $x$.

Two meet-tree $P$, $P'$ are {\it isomorphic} if they are isomorphic as
posets; in particular, an isomorphism $f$ from $P$ to $P'$ preserves
meets, that is $f(x\wedge y)=f(x)\wedge f(y)$ for all $x,y\in P$.  A
{\it positive real-valued meet-tree} , valued meet-tree for short, is
a pair $(P, v)$ where $P$ is a meet-tree and $v$ a map from $P$ to
$\R^*$. Two valued meet-trees $(P, v)$, $(P,v')$ are {\it isomorphic}
if there is an isomorphism $f$ from $P$ onto $P'$ such that $v'\circ
f= v$.  A {\it subtree} of a meet-tree $P$ is a subset $P'$ of $P$
such that the meet of two arbitrary elements of $P'$ belongs to $P'$;
a {\it valued subtree} of a valued meet-tree $(P, v)$ is a pair
$(P',v')$ where $P'$ is a subtree and $v':= v_{\restriction P'}$. The
{\it age } of a valued meet-tree $(P, v)$ is the collection of finite
valued meet-trees which are isomorphic to some valued subtree of $P$.

Let $\M=(M,d)$ be a metric space, $r\in R_+$ and $a\in M$,
the {\it closed} {\it ball of center $a$, radius $s$} is the set
$\mathcal{B}_a(s):=\{x\in M \mid d(a,x)\leq s\}$. The {\it diameter}
of a subset $B$ of $E$ is $\delta(B): =sup\{d(x,y): x,y \in B\}$. We
denote by $\mathcal Ball(\M)$ be the collection of closed balls of $M$
and by $Nerv (\M):= \{\mathcal{B}_a(s): a\in M, s\in Spec(\M,a)\}$.

\noindent Notice that $\delta(B_a(s))=s$ whenever $s\in
Spec(\mathbb{M}, a)$, but more genarally let us recall the following
fact.

\begin{lem}\label{trivial}   
If $M$ is an ultrametric space then for every $B\in \mathcal Ball(\M)$
and $a\in B$, $B=\mathcal{B}_a(s)$ where $s:= \delta (B)$,
\end{lem}

We give below a description of ultrametric spaces in terms of valued
trees. A very close description is given by Lemin \cite{lemin2} (who
instead of $Nerv(\M)$ considered $Ball(\M)$).

\begin{thm} \label{characterisation1}
\begin{enumerate}
\item Let $\M:= (M, d)$ be an ultrametric space, then the pair $(P,
v)$, where $P: =(Nerv(\M), \supseteq)$, $v$ is the diameter function,
is a valued ramified meet-tree such that every element is below some
maximal element and  the map $\delta: Nerv(\M) \rightarrow Spec(\M)$ is
strictly decreasing, $\delta (X)=0$ for every $X\in M':= \max (P)$ and
$d(x,y)= \delta (\{x\}\wedge \{y\})$ for every $x,y\in M$.

\item Conversely, let $(P, v)$a valued ramified meet-tree such that
every element is below some maximal element of $P$ and the map $v:
P\rightarrow \R_+$ is strictly decreasing with $v (x):= 0$ for each
maximal element of $P$. Then the map $d$ defined on $M':= \max (P)$ by
$d(x,y): = v (x \wedge y)$ is an ultrametric distance and $Nerv (\M')=
up(P)_{\restriction M'}$ where $up(P)_{\restriction M'}:=\{M'\cap
\uparrow x : x\in P\}$.
\item  The two correspondences are inverse of each other. 
\end{enumerate}

\end{thm}

\begin{proof} 
1) According to Lemma~\ref{trivial}, balls are disjoint or comparable
w.r.t. inclusion, hence $P$ is a tree. Since $\{x\}\in P$ for every
$x\in M$, $P$ is ramified and every element is below some maximal
element. Let $B ,B'\in P$. Pick $a\in B$, $a'\in B'$ and set $r:=
d(a,a')$. It is easy to se that $\mathcal{B}_a(r)= B\wedge B'$, hence
$P$ is a meet-tree. The properties of $\delta$ follows from Lemma~\ref
{trivial}.

2) a) $d$ is an ultrametric distance~: Let $x\in M'$. We have
$d(x,x):=v(x\wedge x)=v(x)=0$.  If $x\not =y$ then, since $v$ is
strictly decreasing, $d(x,y):= v(x\wedge y)>v(x)= 0$.  Clearly
$d(x,y)= d(y, x)$.  Let $x,y,z\in M'$.  Since $P$ is a tree, $x\wedge
z$ and $y\wedge z$ are comparable. Suppose $x\wedge z\leq y\wedge z $.
Then $x\wedge z\leq x\wedge y $. Since $v$ is decreasing, we have
$d(x,y)\leq d(x,z) \leq \max  \{d(x,z), d(y, z)\}$.

b) $Nerv (\M')= up(P)_{\restriction M'}$~:

Let $B:=M'\cap \uparrow x \in up(P)_{\restriction M'}$, $r:= v(x)$ and
$y\in B$.

{\bf Claim 1.}  $B=\mathcal{B}_y(r)$ and $r\in Spec(\M', y)$.  Thus $B\in Nerv(\M')$. 

Indeed, let $z\in B(y, r)$, that is $v(y \wedge z)\leq r$.  Since
$x\leq y$ and $y\wedge z\leq y$, $x$ and $y\wedge z$ are comparable,
since $v$ is strictly decreasing $x\leq y\wedge z$ hence $z\in
B$. Conversely, if $z\in B$ then $x\leq y\wedge z $ thus, since $v$ is
strictly decreasing, $d(y,z):= v(y\wedge z)\leq v(x)= r$ proving $z\in
\mathcal{B}_y(r)$. Thus $B=\mathcal{B}_y(r)$ as claimed. Since $P$ is
ramified and every element of $P$ is below some element of $M'$, there
is some $z\in M'$ such that $x=y\wedge z$. Clearly, $z\in B$ and
$r=d(y, z)$ thus $r\in Spec(\M', y)$.
 
Let $B: = B(y, r)\in Nerv(\M')$ with $r\in Spec(\M', y)$
 
{\bf Claim 2. }  $B\in  up(P)_{\restriction M'}$. 
 
Indeed, since $r\in Spec(\M', y)$ there is some $z\in M'$ such that
$d(y,z)= r$. Let $x:= y\wedge z$. Since $v(x)= r$ we get $B= \uparrow
x\cap M'\in up(P)_{\restriction M'}$ from the previous claim.

3) We simply note that if $P:= (Nerv(\M), \supseteq)$ then, for $M':=
\max (P)$, $P$ is isomorphic to $(up(P)_{\restriction M'}, \supseteq)$;
moreover, if $v: P\rightarrow \R_+$ is the diameter function
associated to $\M$, then $v(x)= \delta' (M'\cap \uparrow x )$ where
$\delta$ is the diameter function associated to the metric defined on
$M'$ in part 2.\end{proof}

\begin{lem} \label{ages1}
Two ultrametric spaces have the same age if and only if the
corresponding valued trees have the same age.
\end{lem}

The verification is immediate.

The {\it reduced valued tree associated to} an ultrametric space $\M$
is the pair $(P', v')$ where $P':= P\setminus \max (P)$ and $v':=
v_{\restriction P'}$. The age of the reduced valued tree does
not determine the age of the tree, because the information about the
degree, in $P$, of terminal nodes in $P'$ is missing. With this
information added, we have easily:

\begin{lem} \label{ages2}
If two reduced valued trees are isomorphic via a map which preserves
the degree of the original trees then the ultrametric spaces have the
same age.
\end{lem}


Let $\lambda$ be a chain and let $\overline a:= (a_{\mu})_{\mu\in
\lambda }$ such that $2\leq a_{\mu} \leq \omega$. Set
$\omega^{[\overline a]}:=\{\overline b:= (b_{\mu})_{\mu\in \lambda} :
\mu\in \lambda \Rightarrow b_{\mu} <a_{\mu}$ and $supp(\overline b):=
\{\mu <\omega: b_{\mu} \not = 0\}$ is finite $\}$.  If
$a_{\mu}=\omega$ for every $\mu \in \lambda$, the set
$\omega^{[\overline a]}$ is usually denoted $\omega^{[\lambda]}$. Add
a largest element, denoted $\infty$ to $\lambda$. Given $\overline b,
\overline c \in \omega^{[\overline a]}$, set $\Delta(\overline b,
\overline c):=\infty$ if $\overline b= \overline c$, otherwise
$\Delta(\overline b, \overline c):=\mu$ where $\mu$ is the least
member of $\lambda$ such that $ b_{\mu}\not = c_{\mu}$.

Suppose $\lambda$ be countable. Let $w: \lambda\cup
\{\infty\}\rightarrow \R_+$ be a strictly decreasing map such that
$w(\infty)=0$, let $d_w: = w\circ \Delta $ and let $V$ be the image
of $w$.  For $\mu \in  \lambda \cup \{\infty\}$ set
$\downarrow^*\mu:=\downarrow \mu\setminus \{\mu\}$. Let 
$P':=\{f_{\restriction\downarrow^*\mu}: f\in \omega^{[\overline a]},  \mu\in \lambda \cup \{\infty\}\}$ ordered by extension and let $v'(f_{\restriction\downarrow^*\mu}):= w(\mu)$. 
  
 We have the following property, which is easy to check.
 
 \begin{lem} \label{homogeneous0}
The pair $\M:= (\omega^{[\overline a]}, d_w)$ is an ultrametric space,
$ Spec(\M)= V$ and the  valued tree associated
to $\M$ is isomorphic to $(P', v')$.  
\end{lem}

We say that $\M$ is {\it point-homogeneous} if the automorphism group
of $\M$ acts transitively on $\M$.

\begin{thm}\label {homogeneous1}
Let $\M$ be a countable ultrametric space, $P:= (Nerv(\M),
\supseteq)$, $v: P\rightarrow \R_+$ where $v(B):= \delta(B)$, $M':=
\max (P)$.  The following properties are equivalent:

\begin{enumerate}[{(i)}]
\item  $\M$ is isometric to some $(\omega^{[\overline a]}, d_w)$ .
\item $\M$  is homogeneous;
\item $\M$ is point-homogeneous;
\item 
\begin{enumerate}
\item  $v(x)=v(y)\Rightarrow d_P(x)=d_P(y)$ for every $x,y\in P$;
\item  $ v[\downarrow x]= v[\downarrow y]$ for every $x, y\in \max (P)$.
\end{enumerate}
\end{enumerate}
\end{thm}

\begin{proof}
$(i)\Rightarrow (iv)$ Let $\M:= (\omega^{[\overline a]}, d_w)$. 
 According to Lemma~\ref{homogeneous0}, the valued tree associated to $\M$ is isomorphic to $(P', v')$.  Condition $(b)(iv)$ immediately follows. Let $x:=
f_{\restriction\downarrow^*\mu}\in P'$;  if $\mu= \infty$ then $d_{P'}(x)= 0$, otherwise $d_{P'}(x)=a(\mu)$. Thus
Condition $(a)(iv)$ holds too.

$(ii)\Rightarrow (iii)$ Trivial

$(iii)\Rightarrow (iv)$ Suppose $\M$ point homogeneous. First,
Condition $(b)(iv)$ holds. Indeed, let $x,y\in M':=\max  (P)$.  Then
$x:= \{x'\}$ and $y:= \{y'\}$, with $x',y'\in \M$. Let $f$ be an
isometry from $\M$ onto itself such that $f(x')=y'$. Then $Spec(x',
\M)= Spec(y', \M)$ and the result follows. Next, Condition $(a)(iv)$
holds. Let $x:= B\in P, y:= C\in P$ and $r:= v(x)=v(y)$.  Pick $x'\in
B$, $y'\in C$.  Let $f$ be an isometry from $\M$ onto itself such that
$f(x')=y'$. Then $f(B)=C$. For two elements $x',y'$ of $B$, set
$x'\equiv y'$ if $d(x',y')<r$. This relation is an equivalence
relation.  whose number of classes is the degree of $x:=B$ in the
poset $P:=Nerv(\M)$. The desired conclusion follows. 

 $(iv)\Rightarrow
(ii)$ Let $f$ be an isometry from a finite subset $A$ of $M$ onto a
subset $B$ of $M$. Let $x\in M\setminus A$. We prove that $f$ extends
to an isometry defined on $A\cup \{x\}$.  If $A$ is empty, we may send
$x$ onto any element $b$ of $M$. If $A$ is non-empty, set $r: =
\min (\{d(x,y): y\in A\})$. In order to extend $f$ we only need to send
$x$ onto some $b\in M$ such that $f(B(x, r))= B(b, r)\cap f(A)$.
There is some $u\in P$ such that $x\wedge x'=u$ for all $x' \in B(x,
r)\cap A$ and moreover $v(u)=r$. Select $y\in f(B(x, r))$. Since
$v[\downarrow x]= v[\downarrow y]$ there is some $u'\in \downarrow y$
such that $v(u')= r$. Since $d_P(u)= d_P (u')$, there is $b\in M$ such
that $y'\wedge b= u'$ for all $y'\in f(A)$. Such an element will do.

$(ii)\Rightarrow (i)$. Let $\lambda:= Spec(\M)\setminus \{0\}$ ordered
with the dual of the order induced by the natural order on $\R$, let
$w: \lambda\cup\{\infty \}\rightarrow \R_+$ with $w(x):=x$ for $x\in \lambda$ and $w(\infty):= 0$ and let $\overline a: \lambda\rightarrow
\omega+ 1$ such that $\overline a\circ w= d_P$ (such a map exists
because of (iv) Condition 1).  

{\bf Claim} $\M$ is isometric to
$(\omega^{[\overline a]}, d_w)$.
According to the implications $(i)\Rightarrow (iv)\Rightarrow (ii)$ proved  above,  $(\omega^{[\overline a]}, d_w)$ is
homogeneous.  Since $\M$ is homogeneous, it suffices to prove that
$(\omega^{[\overline a]}, d_w)$ and $\M$ have the same age to get the
desired conclusion. From the implication $(iii)\Rightarrow (iv)$, the reduced valued trees associated to
$(\omega^{[\overline a]}, d_w)$ and $\M$ are isomorphic by an
isomorphism which preserves the degree. From Lemma~\ref {ages2},
$(\omega^{[\overline a]}, d_w)$ and $\M$ have the same age.

\end{proof}

\begin{proposition}\label {homogeneous2}
The space $(\omega^{[\lambda]}, d_w)$ is the countable homogeneous
ultrametric space $\mathbb Ult_V$ associated with $V$.
\end{proposition}

\begin{proof} 
We only neeed to prove that every finite ultrametric space $\mathbb
M:=(M,d)$ with spectrum included into $V$ embeds isometrically into
$(\omega^{[\lambda]}, d_w)$. We argue by induction on the number $n$
of elements of $M$. If $n\leq 1$, the result is obvious. Suppose
$n\geq 2$.  Let $x\in M$. We may suppose that there is an isometric
embedding $f$ of $\mathbb M_{-x}:= \mathbb M_{\restriction {M\setminus
\{x\}}}$ into $(\omega^{[\lambda]}, d_w)$. We prove that $f$ extends
to $\mathbb M$.  Set $r: = \min (\{d(x,y): y\in M\setminus \{x\}\})$ and
$\mu\in \lambda$ such that $w(\mu)= r$. In order to extend $f$ we only
need to find some element $b\in \omega^{[\lambda]}$ such that $f(B(x,
r))= B(b, r)\cap f(M\setminus \{x\})$. For every $\overline b',
\overline b''\in f(B(x, r))$ we have $b'_{\mu'}=b''_{\mu'}$ for all
$\mu'<\mu$. Select $b\in \omega^{[\lambda]}$ such that $b_{\mu'}= b'_{
\mu'}$ for all $\mu'<\mu$ and $b_{\mu}\in \omega\setminus\{b'_{\mu}:
\overline b'\in f(B(x, r)\}$. \end {proof}

\subsection{Indivisible ultrametric spaces}

\begin{defin}\label{defin:ring} Let $\mathbb{M}:=(M;d)$ be a metric space,   
$a\in M$   and $0\leq r<s $.  Then 
\[
\mathcal{R}_a(r,s):=\{x\in M \mid r\leq d(a,x)< s\}.
\]
\end{defin}

\begin{lem}\label{lem:spec}
Let $\mathbb{M}=(M,d)$ be an indivisible ultrametric space. Then the
spectrum of every element of $M$ is dually well founded.
\end{lem}

\begin{proof}
Let $a\in M$. Suppose for a contradiction that
$r_0=0<r_1<r_2<r_3<\dots$ is an infinite sequence of reals in the
spectrum of $a$. Let $s$ be its supremum. Cover $M$ by a family
$\mathcal B:= \{\mathcal{R}_{a_{\alpha}}(0,s) : \alpha< \kappa\}$ of
open balls of radius $s$ such that $a_{\alpha}\not \in M_{\alpha}:=
\cup\{\mathcal{R}_{a_{\beta}}(0,s): \beta <\alpha\}$ (with the
convention that if $s=\infty$ then $\mathcal B$ consists of
$M$). Since $d$ is an ultrametric distance, these balls are pairwise
disjoint and therefore, the rings
$\mathcal{R}_{a_{\alpha}}(r_{i},r_{i+1})$ make-up a partition of
$M$. Let:
\[
\mathcal{E}:=\bigcup_{\alpha<\kappa, i\in \omega}\mathcal{R}_{a_{\alpha}}(r_{2i},r_{2i+1}) \text{  and  }
\mathcal{O}:=\bigcup_{\alpha<\kappa, i\in \omega}\mathcal{R}_{a_{\alpha}}(r_{2i+1}, r_{2i+2})
\]
and let $f$ be an isometry of $M$ into $M$. Let $\alpha<\kappa$ and
$i\in\omega$ so that $f(a)\in
\mathcal{R}_{a_\alpha}(r_i,r_{i+1})$. Let $b\in M$ with
$d(a,b)=r_{i+1}$.

Then $d(f(a),f(b))=r_{i+1}$ and because $d(f(a),a_\alpha)<r_{i+1}$ it
follows that that $d(a_\alpha,f(b))=r_{i+1}<s$.  Thus $f(b)\in
\mathcal{R}_{a_\alpha}(r_{i+1},r_{i+2})$.
\end{proof}

\begin{cor}\label{lem:specmax} If an ultrametric space is indivisible then  the collection of balls, once ordered by inclusion,  is dually well-founded and the diameter is  attained. \end{cor}

\begin{proof} Let $(B_n)_{n<\omega}$ be an
increasing sequence of balls of an ultrametric space $\M:= (M,d)$. Pick $a\in \cap \{B_n: n\in \}$. Since $\M$  is ultrametric, $a$ is the center of each $B_n$ thus their  radii belong to the spectrum of $a$.  If $\M$ is indivisible,  then from Lemma~\ref {lem:spec} above  $Spec(\M, a)$ is dually well-founded, thus the sequence is eventually constant. Let   $s$ be the maximum of $Spec(\M, a)$. Let $x,y\in M$. We have $d(x,y)\leq \max (\{d(x,a), d(y,a)\})\leq s$, hence $s$ is the maximum of the spectrum of $\mathbb{M}$, that is the diameter of $\M$. \end{proof}



%

\begin{thm}\label{thm: ultrahomind}
Let $\M$ be a denumerable  ultrametric space.
The following properties are equivalent:
\begin{enumerate}[{(i)}]
\item $\M$ is isometric to some $\mathbb Ult_V$, where $V$ is dually  well-ordered; 
\item $\M$ is point-homogeneous, $P:=(Nerv(\M), \supseteq )$ is well founded and the degree of every non maximal element is infinite;
\item $\M$ is  homogeneous and indivisible; 
 
\end{enumerate}
\end{thm}

\begin{proof}
$(i)\Rightarrow (ii)$ By definition, $\mathbb Ult_V$ is homogeneous, hence point-homogeneous. In fact, according to Proposition~\ref {homogeneous2},  $\M$ is isometric to  some $(\omega^{[\lambda]}, d_w)$ where $\lambda$ is a well-ordered chain.  Thus, from Lemma~\ref {homogeneous0}, $P:=(Nerv(\M), \supseteq )$ is well-founded and the degree of every non maximal element is infinite.

$(iii)\Rightarrow (i)$ Suppose that $(iii)$ holds. 
Theorem~\ref{homogeneous1} asserts that $\M$ is isometric to some
$(\omega^{[\overline a]}, d_w)$. Since $\M$ is indivisible, it follows
from Lemma~\ref {lem:spec} that $V:= Spec(\M)$ is well-founded, hence
we may suppose that $\lambda$ is an ordinal.  To conclude it suffices
to prove that $a_{\mu}=\omega$ for every $\mu < \lambda$. Let
$\mu<\lambda$; set $r:= w(\mu)$. First, observe that $M= \cup\mathcal
B$ where $\mathcal B$ is a collection of pairwise disjoint balls, all
of diameter $r$. Next, each member $B$ of
$\mathcal B$ is the union of $a_{\mu}$ balls $B_i$ each of smaller
diameter than $r$.  Indeed, since $\M$ is point-homogeneous, all balls having
the same radius are isometric spaces, thus  it suffices to prove this
property for the ball $B:= \mathcal B_0(r)$, where $0$ is the ordinal
sequence which only takes value $0$. This is easy: set $\overline
x_i:= (b_{\nu})_{\nu<\lambda}$ where $i<a_n$, $b_{\nu}= 0$ if $\nu
\not = \mu$ and $b_{\mu}:= i$ otherwise, set $r^+:= w(\mu^+)$ where
$\mu^+:=\mu+1$ if $\mu+1<\lambda$ and $\mu^+:=\infty$ otherwise, then
$B:= \cup\{B(\overline x_i, r^+): i<a_{\mu}\}$. With these two
observations we have $M=\cup\{M_i: i<a_{\mu}\}$ where $M_i:= \cup \{B_i:
B\in \mathcal B\}$. Clearly, there is no isometry from $\M$ into an
$\M_i$ hence if $a_{\mu}<\omega$, $\M$ cannot be indivisible.

$(ii)\Rightarrow (iii)$ According to Theorem~\ref{homogeneous1}, $\M$ is homogeneous. Let us show that it is indivisible. 
Let $f: M\rightarrow 2$ be a partition of $M$ into two parts. Set
$\mathcal F_0$ be the set of balls $B\in Nerv(\M)$ such that there is some
isometry $\varphi_B$ from $B$ into$B\cap f^{-1}(0)$ and let $M_0:=
\cup \mathcal F_{0}$.

{\bf Claim 1} There is an isometry from $M_0$ to $M_{0}\cap
f^{-1}(0)$.

Indeed, let $\mathcal F'_{0}$ be the subset of $\mathcal F_{0}$ made
of its maximal members (w.r.t. inclusion).  Let $\varphi:=\cup \{\varphi _B: B\in
\mathcal F'_{0}\}$.  Since balls are either disjoint or comparable,
$\varphi$ is a map and, since $P:= (Nerv(\M), \supseteq)$ is  well-founded,
$M_{0}= \cup \mathcal F'_{0}$, hence the domain of $\varphi$ is
$M_{0}$.

For $B$ in $Nerv(\M)$, set $Pred(B):=\max  ( \{B': B' \subset B, B'\in
Nerv(\M)\})$.

{\bf Claim 2} If $B\not \in \mathcal F_{0}$ then $Pred(B)\cap \mathcal
F_{0}$ is finite.

Indeed, suppose not.  Then, since the space is point-homogeneous, all
members of $Pred(B)$ have the same radius and there is an isometry
$\psi$ from $B$ into $B$ which transforms each member of $Pred (B)$ to
a member of $Pred(B)\cap \mathcal F_{0}$.  Let $\varphi:=\cup
\{\varphi _{B'}: B'\in Pred(B)\cap \mathcal F_{0}\}$. Then $\varphi$
is an isometry from $\cup ( Pred(B)\cap \mathcal F_{0})$ into $B\cap
f^{-1}(0)$. Consequently, $\varphi\circ \psi$ is an isometry from $B$
into $B\cap f^{-1}(0)$, thus $B\in \mathcal F_{0}$, a contradiction.
  
Suppose that $M\not \in \mathcal F_{0}$. We construct an isometry $h$
from $M$ into $f^{-1}(1)\setminus M_{0}$ as follows. We start with an
enumeration $(x_n)_{n<\omega}$ of the elements of $M$.  According to
Claim 1, $M\setminus M_{0}\not =\emptyset$. We may also suppose that
it contains an element of $f^{-1}(1)$ (otherwise the union of the
identity map on $M\setminus M_{0}$ and an isometry as constructed in
Claim 1, is an isometry from $M$ into $f^{-1}(0)$). Let $y_{0}$ such
an element. We set $h(x_0):= y_{0}$.

Suppose $h$ be defined for all $m$, $m<n$. Let $p:= \min  (\{d(x_m,
x_n): m<n\})$. Let $I:= \{i, i<n: d(x_i, x_n):=p\}$ . Let $B:=
\mathcal B_{h(i)}( p)$ for $i\in I$. This set does not depend upon the
choice of $i$.  Since $h(i)\in f^{-1}(1)\setminus M_{0}$, $B\not \in
\mathcal F_{0}$. For each $i\in I$ let $B'_{i}$ such that $h(i)\in
B'_{i}\in Pred(B)$. According to Claim 2, there is some $B''\in
Pred(B)\setminus \mathcal F_{0}$ which is distinct from all the
$B_{i}'s$. As in our first step, $B''\setminus M_{0}$ is nonempty and
in fact contains an element, say $y_n$ of $f^{-1}(1)$. We set
$h(x_n):= y_n$.

\end{proof}

\section{Divisibility of  metric spaces}

The sequence $a_0,a_1,\dots, a_{n-1},a_n$ of elements in a metric
space $\M:= (M;d)$ is an {\em $\epsilon$-chain joining $a_0$ and
$a_n$} if $d(a_i,a_{i+1})\leq \epsilon$ for all $i\in n$.  The space
$\M$ is {\em Cantor connected} if any two of its elements can be
joined by an $\epsilon$-chain for any $\epsilon>0$.  The {\it Cantor
connected component of an element $a\in M$} is the largest Cantor
connected subset of $M$ containing $a$. The space $\M$ is {\em totally
Cantor disconnected} if the Cantor connected component of every $a$
reduce to $a$.  See \cite{Lemin} for more details and references.

For $a\in M$ let $\lambda_\epsilon(a)$ be the supremum of all reals
$l\leq 1$ for which there exists an $\epsilon$-chain $a_0,a_1,\dots,
a_{n-1},a_n$ with $d(a_0,a_n)\geq l$ containing $a$. (The condition
$l\leq 1$ saves us from having to consider the special case $\infty$.)
Let
\[
\lambda(a):=\sup\{l\in \mathbb{R} \mid \forall \epsilon>0\,
(\lambda_\epsilon(a)\geq l) \}.
\]
A space $(M;d)$ is {\em restricted} if $\lambda(a)=0$ for all $a\in
M$. It follows that every restricted space is totally Cantor
disconnected. There are totally Cantor disconnected spaces which are
not restricted. Here is an example with a finite diameter~:

\begin{example}\label{ex:1}
Let $(M;d)$ be the metric space so that:
\vskip -20pt

\begin{enumerate}
  \item $M=\{(0,0)\}\cup\{(m,n)\in \N\times\N\ |\ m<n\}$ 
  \item $d((0,0),(m,n))=\frac m n$
  \item $d((m_1,n),(m_2,n))=\frac{|m_1-m_2|}n$
  \item $d((m_1,n_1),(m_2,n_2))=\frac{m_1}{n_1}+\frac{m_2}{n_2}$ when $n_1\neq n_2$.
\end{enumerate}
\end{example}

\begin{lem}\label{lem:2}
Let $c\in M$ and $0\leq r_0<r_1<r_2<r_3 $ and $a\in
\mathcal{R}_c(r_0,r_1)$ and $b\in\mathcal{R}_c(r_2,r_3)$ then:
\begin{enumerate}
\item $d(a,b)>r_2-r_1$.
\item $d(x,y)<2r_2$ for all $x,y\in \mathcal{R}_c(r_1,r_2)$.
\item If $0<\epsilon<\min\{r_1-r_0,r_3-r_2\}$ and
$x_0,x_1,x_2,\dots,x_{n-1}$ is an $\epsilon$-sequence with $x_i\not\in
\mathcal R_c(r_0,r_1)\cup \mathcal R_c(r_2,r_3)$ for all $i\in n$ but with $x_i\in
\mathcal{R}_c(r_1,r_2) $ for at least one $i\in n$, then $x_i\in
\mathcal R_c(r_1,r_2)$ for all $i\in n$.
\item Let $f$ be an isometry of $M$ with $f[M]\cap
(\mathcal{R}_c(r_0,r_1)\cup \mathcal{R}_c(r_2,r_3))=\emptyset$ and let
$z\in M$ with $\lambda(z)> 2r_2$. Then $f(z)\not\in
\mathcal{R}_c(r_1,r_2)$.
\end{enumerate}
\end{lem}
\begin{proof}
Items 1 and 2 follow from the triangle inequality. Item 3 follows from
item 1 and item 4 follows from items 2 and 3.
\end{proof}

\begin{defin}\label{defin:stripes}
Let $c\in M$ and $0<l$. Then 
\[
\mathcal{E}_c(l):=\bigcup_{ n \geq 2, \; n \; \mathrm{ even}}
\mathcal{R}_c\bigl(\frac{l(n-1)}{n}, \frac{ln}{n+1}\bigr)
\]
and
\[
 \mathcal{O}_c(l):=\bigcup_{n\;   \mathrm{ odd}}
 \mathcal{R}_c\bigl(\frac{l(n-1)}{n}, \frac{ln}{n+1}\bigr). 
 \]
\end{defin}

\begin{thm}\label{lem:basic_divisibility}
Let $\mathbb{M}=(M;d)$ be a countable metric space. If there exists an
element $a\in M$ with $\lambda(a)>0$ then $\mathbb{M}$ is divisible.
\end{thm}

\begin{proof}
Since $M$ is countable, it can be covered by a family of pairwise
disjoint open balls with radius less that $\frac {\lambda(a)} {2}$.
In fact, there exists a subset $C$ of $M$ and for every $c\in C$ a
positive real $l_c$ so that:

\begin{enumerate}
\item $l_c\not=d(x,y)$  for every $c\in C$ and $x,y\in M$.
\item $2l_c<\lambda(a)$ for every $c\in C$.
\item For every element $x\in M$ there is one and only one element
$c\in C$ with $x\in \mathcal{R}_c(0,l_c)$.
\end{enumerate}

(After enumerating $M$ into an $\omega$-sequence
$m_0,m_1,m_2,m_3,\dots$ such a set $C$ and function $l$ can be
constructed step by step exhausting all of the elements of $M$.)

Let 
\[
\mathcal{E}:= \bigcup_{c\in C} \mathcal{E}_c(l_c)  \text{  and   } 
\mathcal{O}:=\bigcup_{c\in C}\mathcal{O}_c(l_c).
\]
Then $\mathcal{E}\cup \mathcal{O}=M$ and $\mathcal{E}\cap \mathcal{O}=\emptyset$.

Assume for a contradiction that there is an isometry $f$ which maps
$M$ into $\mathcal{E}$. Then there is a $c\in C$ so that $f(a)\in
\mathcal{E}_c(l_c)$. But this is not possible according to Lemma
\ref{lem:2} item 4. Similarly it is not possible that $f$ maps $M$ into
$\mathcal{O}$.
\end{proof}

\begin{cor}\label{cor:basic_divisibility}
A countable metric space which is indivisible is restricted and hence
totally Cantor disconnected.
\end{cor}

The second part of the conclusion of the corollary above extends to
uncountable metric spaces.

\begin{thm} \label{uncountable}
Let $\M$ be a metric space and $r$ be a positive real, then there is a
partition into two parts $A_{0}$ and $A_{1}$ which contains no Cantor
connected subspace of diameter larger than $r$.
\end{thm}

\begin{lem}
Le $\M$ be a metric space and $r$ be a positive real number. Then
there is a sequence $(E_\mu)_ {\mu< \lambda}$ such that:

\begin{enumerate}

\item $E_0 = \emptyset$ and each $E_\mu$ is open in $M$

\item the sequence is strictly increasing and continuous, that is
$E_\mu$ is the union of $E_{\nu} \text{ for } \nu<\mu$ if $\mu$
is a limit ordinal,

\item the union covers $M$

\item $F_\mu:= E_{\mu+1}\setminus {E_\mu}$ has diameter at most $r$  and decomposes into two sets $A_{\mu,
0}$ and $A_{\mu, 1}$ such that each Cantor connected subspace  $Y$
of $A_{\mu, i}$ is contained into some subset $\mathcal B_Y$ of $F_\mu$ such that $d(y, A_{\mu, i}\cap \mathcal B_Y)\geq \epsilon_Y$ for some $\epsilon_Y>0$ and
every $y\in M\setminus (\mathcal B_Y \cup E_\mu)$.

\end{enumerate}
\end{lem}

\begin{proof}
Suppose the sequence defined for all $\nu$, $\nu<\mu$. If $\mu$ is a limit
ordinal, set $E_\mu:= \bigcup \{E_\nu : \nu<\mu\}$. If $\mu$ is a successor,
say $\mu:= \nu+1$, pick $x\in E':= M\setminus E_\nu$, set $\mathcal R'_x( 0,
r/2):= \{y\in E': d(x,y)< r/2\}$ and set $E_{\mu}:= E_\nu\cup \mathcal
\mathcal R'_x( 0, r/2)$. Decompose $\mathcal R'_x( 0, r/2)$ into
countably many crowns $\mathcal{R'}_x\bigl(\frac{r(n-1)/2}{n},
\frac{rn/2}{n+1}\bigr)$ as in the proof of Theorem~\ref
{lem:basic_divisibility} , the union of the even ones gives $A_{\nu,0}$,
the rest gives $A_{\nu,1}$.

\noindent Finally any Cantor connected subspace  $Y$ of $A_{\mu, i}$ must be
included into  $\mathcal{R'}_x\bigl(\frac{r(n-1)/2}{n},
\frac{rn/2}{n+1}\bigr)$  for some $n$, and therefore the required $\mathcal B_Y$
may be taken to be $\mathcal{R'}_x\bigl(0,
s_Y \bigr)$ with $\frac {rn/2}{n+1}<s_Y< \frac {r(n+1)/2}{n+2}$ with  $\epsilon_Y:= s_Y-\frac {rn/2}{n+1}$.

\end{proof}

\noindent {\bf Proof of  Theorem~\ref{uncountable}}
Let $A_{i}:= \bigcup \{A_{\mu, i}: \mu<\lambda\}$.
Then $A_{i}$ contains no  Cantor connected subspace $X$ of diameter larger than $r$.

Indeed, suppose the contrary.  Let $\mu$ be minimum such that $E_{\mu}$
meets  $X$. Clearly $\mu$ is a successor, say $\mu=\nu+1$.  Let $x\in X_\nu:= X\cap
F_\nu$. Let $Y$ be the Cantor connected component of $x$ in  $A_{\nu, i}$ and let $\mathcal B_Y$ given by the above lemma.
{\bf Claim}  $X\subseteq \mathcal B_Y$. Indeed suppose not, let
$y\in X\setminus \mathcal B_Y$, let $\epsilon$,  $0<\epsilon< \epsilon_Y$ and  $x_0:=x,\dots, x_k, \dots, x_n=y$ be an
$\epsilon$ path contained in $X$.  Let $\ell$ be least index such that $x_\ell\not \in
\mathcal B_Y$. From $x_{\ell-1} \in
A_{\mu, i}\cap \mathcal B_Y$, we get  $d(x_\ell, A_{\mu, i}\cap \mathcal B_Y)< \epsilon_Y$. A
contradiction.

Since $X\subseteq B_Y\subseteq F_\mu$, the  diameter of $X$ is at most $r$. The proof is complete.

\begin{defin}\label{defin:semi}
Let $\mathbb{M}:=(M;d)$ be totally Cantor disconnected. Then 
\[
d^\ast(x,y):=\inf\{\epsilon>0 \mid \text{ there exists an
$\epsilon$-sequence containing $x$ and $y$}\}.
\]
\end{defin}

\begin{lem}\label{lem:semi}
Let $\mathbb{M}:=(M;d)$ be totally Cantor disconnected. Then
$\mathbb{M}^\ast:=(M;d^\ast)$ is an ultrametric space.
\end{lem}

\begin{proof}
Let $x,y,z\in M$ with $d^\ast(x,y)\geq d^\ast(x,z)\geq
d^\ast(y,z)$. Then for every $\epsilon>d^\ast(x,z)$
there are $\epsilon$-sequences 
joining $x$ to $z$ and $z$ to $y$, then one
joining $x$ to $y$, hence $d^\ast(x,y)\leq\epsilon$.
Thus $d^\ast(x,y)\leq d^\ast(x,z)$.
\end{proof}

(See \cite{Lemin}, Theorem 1 and Lemma 8.)

\begin{thm}\label{thm:ultra}  
Let $\mathbb{M}:=(M;d)$ be a countable homogeneous indivisible metric
space then $\mathbb{M}^\ast$ is an homogeneous indivisible ultrametric
space.
\end{thm}

\begin{proof} 
Since $\M$ is indivisible it is totally Cantor disconnected, hence
$d^\ast$ is well defined. Since $\M$ is homogeneous then
$d(x,y)=d(x',y')$ implies $d^\ast (x,y)= d^\ast (x',y')$ for all
$x,y,x',y'\in M$.  From this property, every local isometry of $\M$ is
a local isometry of $\M^\ast$. Hence, since $M$ is indivisible,
$M^\ast $ is indivisible. Since every automorphism of $M$ is an
automorphism of $M^\ast$, $M^\ast $ is point-homogeneous. According to
Theorem~\ref{homogeneous1}, $M^\ast$ is homogeneous.
\end{proof}

\begin{thm}
Let $\mathbb{M}$ be a homogeneous metric space and $V:=
Spec(\mathbb{M})$.  If $\mathbb{M} $ is totally Cantor disconnected
and every three element metric space $\mathbb T$ with $Spec(\mathbb
T)\subseteq V$ embeds into $\mathbb{M}$ then the set $V\setminus
\{0\}$ is either contained into an interval of the form $[a\rightarrow
+\infty)$ for some $a\in \R_+\setminus \{0\}$ or into an union of
intervals of the form $\cup \{ [a_{2(n+1)}, a_{2n+1}]: n<\omega\} \cup
[a_0\rightarrow +\infty)$ where $\{a_n : n<\omega\}$ is a sequence
such that $a_{2n+1}\leq \frac {a_{2n}}{2}$.
\end{thm}

\begin{proof}
{\bf Claim} For every $w\in V^\ast := Spec(\M^\ast )$, $]\frac{w}{2}, w [\cap V=\emptyset$.

Suppose the contrary.  Pick $r\in ]\frac{w}{2}, w [\cap
V=\emptyset$. Since $w\in V^\ast $, we may find $x,y$ such that
$d^\ast(x,y)= w$. Let $n<\omega$ and $\epsilon:= 2r$, then there is an
$\epsilon$-sequence $x_0,\dots, x_n$ containing $x, y$. For $i<n$, let
$\mathbb T_i := (\{x_i, x_{i+1}, z_i\}, d_i) $ where $d_i(x_i,
x_{i+1}):= d(x_i, x_{i+1})$, $d_i(x_i,z_i)=d_i(x_{i+1}):= r$. Each
$\mathbb T_i$ is a metric space whith spectrum included into $V$,
hence can be isometrically embedded into $\mathbb{M}$. Since $M$ is
homegenous, we may suppose that $z_i\in M$ and that the embedding is
the inclusion. By adding the $z_i'$'s to the $x_i's$ we get a
$r$-sequence containing $x$ and $y$. Since $r<w$ this gives a
contradiction.

Since every element of $V^\ast $ is the infimum of elements of $V$ it
also follows that $]\frac{w}{2}, w [\cap V^\ast =\emptyset$.

Let $\alpha:= Inf (V\setminus \{0\})$. If $\alpha \not=0$ set $a:=
\alpha$; in this case $V\setminus \{0\}\subseteq [a\rightarrow
+\infty)$. If $\alpha=0$ then, since every element de $V\setminus
\{0\}$ majorizes some element of $V^\ast \setminus \{0\}$ it follows
that $Inf (V^\ast \setminus \{0\})=0$ too.  Let $\{a_{2n} :
n<\omega\}$ be a strictly decreasing sequence of elements of $V^\ast$
which converges to $0$.  Set $a_{2n+1}:= a_{2n}$. From the Claim
$]\frac{a_{2n}} {2}\; a_{2n}[\cap V=\emptyset$, hence $a_{2n+2}\leq
a_{2n+1}$. The rest follows.  

\end{proof}

\begin{thm}\label{thm:unbounded} 
Every unbounded metric space is divisible.
\end{thm}

\begin{proof}
Let $\mathbb{M}:=(M;d)$ be an unbounded metric space.  Construct a
sequence of reals $r_0, r_1, r_2, \dots $ and a sequence $a_0, a_1,
a_2, \dots $ of elements of $M$ so that for every integer $i\in\\N$

\begin{enumerate}
\item $d(a_0, a_{i+1})>2r_i$.
\item $d(a_0,a_{i+1})+r_i< r_{i+1}$.
\end{enumerate}

Let $r_0:=0$ and $a_0\in M$ be arbitrary. Suppose that $(r_i :
i\leq n)$ and $(a_i: i\leq n)$ have already been constructed.  From
the fact that $\mathbb{M}$ is unbounded, we can find $a_{n+1}\in M$
such that $d(a_0,a_{n+1})>2r_n$.  Next, choose $r_{n+1}> d(a_0,a_{n+1})+
r_n$. Note that the set $\{r_i : i\in \N\}$ such constructed
is unbounded.

Let, given any $c\in M$,
\[
\mathcal{E}:=\bigcup_{i\in \N}\mathcal{R}_{c}(r_{2i},r_{2i+1}) \text{ and  }
\mathcal{O}:=\bigcup_{i\in \N}\mathcal{R}_{c}(r_{2i+1}, r_{2i+2}).
\]
We prove that there is no isometric embedding of $M$ into
$\mathcal{E}$ or into $\mathcal{O}$.

Let $f$ be an isometric embedding of $M$ into $M$. Let $i$ be minimal
so that $d(c, f(a_0))<r_i$;
notice that $i>0$ and $f(a_0)\in
\mathcal{R}_{c}(r_{i-1}, r_i)$. 
We have:
\begin{align*}
&d(c,f(a_{i+1}))\geq d(f(a_0),f(a_{i+1}))-d(c,f(a_0))=\\
&=d(a_0,a_{i+1})-d(c,f(a_0))>2r_i-r_i=r_i.
\end{align*}
Also: 
\[
d(c,f(a_{i+1}))\leq d(c,f(a_0))+d(f(a_0),f(a_{i+1}))\leq r_i+d(a_0,a_{i+1})<r_{i+1}.
\]
It follows that $f(a_{i+1})\in \mathcal{R}_{c}(r_i,r_{i+1})$.
Therefore $f[M]$ intersects both $ \mathcal{E}$ and $\mathcal {O}$.

\end{proof}

\section{Divisibility of the bounded Urysohn space}

In \cite{Hjorth}, Hjorth shows that the Urysohn space
$\mathbb{U}_{\Q}$ is divisible, and asks whether the
corresponding bounded space has the same property. We show that it
does, and in fact this generalizes to   bounded Urysohn spaces for which the spectrum $V$ satisfies a density condition. In the sequel $V$ will denote a countable subset of $\R_+$ containing $0$ and satisfying the four-values condition.

\begin{proposition}\label {cantorconnected} 
If for some
$r>0$, $V\cap [0, r]$ is dense in $[0, r]$ then the diameter of each
Cantor connected component of $\U_V$ is at least $r$.
\end{proposition}

\begin{proof}
Let $a\in U_V$ and $\ell\in V\cap (0, r]$.  Let $b \in U$ such that
$d(a,b)=\ell $. For any $n \in \omega$ choose successively $a_0,\dots,
a_ n$ such that:

$a_0:=\frac{\ell} {n}$, $a_{k}:=\frac {\ell} {n}k+\epsilon_{k}\in
V\cap [\frac {\ell} {n}(k-1), \frac {a_{k-1}}{2}]$ for $1\leq k<n$ and
$a_n:= \ell$.

Let $x_0=a$, $x_n=b$, $x_1,x_2,\dots,x_{n-1}$ be elements not in $U$
and $X:=\{x_0,x_1,x_2,\dots,x_n\}$. Let $d^\prime: X\times X
\rightarrow V$ defined by $d^\prime(x_i,x_{i+k})=a_{k}$ for $1\leq
k\leq n$ and $d^\prime(x_i,x_i):= 0$.  With our choice,
$\epsilon_{i+j}\leq \epsilon_i+\epsilon_j$ for all $i,j$ such that
$i+j\leq n$, thus $\X:= (X, d^\prime)\in \mathcal M_V$.

\noindent Hence we can use the mapping extension property of the
Urysohn space and obtain an embedding $f$ of the space $\X$ into $U$
which is the identity map on $x_0$ and $x_n$.

\noindent Since this can be done for any $n$ we conclude that the
Cantor connected component of $a$ contains $b$, hence its diameter is
at least $\ell$. Since this holds for every $\ell\in V\cap (0, r]$,
this diameter is at least $r$.
\end{proof}

\begin{remark} 
If $V$ is residuated,  Proposition~\ref {cantorconnected} follows from Proposition~\ref{fundamental}. Indeed, in this case  $\V:= (V, d_V)$ is a metric space, thus for every $a\in \U_V$ there is  isometric embedding of  $\V$ into $\U_V$ which maps $0$ onto $a$ (Point~\ref{universality} of Section~\ref{relational} together with the homogeneity of $\U_V$). The  Cantor connected component of $0$ contains $V\cap [0,r]$, thus its image has diameter at least $r$. 
Notice that when $V$ is an initial segment of $\Q$, $d_V$
is just its usual distance, in which case Proposition~\ref{fundamental} is not required.
\end{remark}

\begin{cor}\label{thm:Cant_Ury}
The countable homogeneous metric space $\mathbb{U}_{\Q,\leq
1}:=(U;d)$ having all rational numbers less than or equal to one as
distances is Cantor connected.
\end{cor}

\begin{thm}  If for some $r>0$, $V\cap [0, r]$ is dense in $[0, r]$ then  $\U_V$ is divisible. 
\end{thm}

\begin{proof}
Follows from Proposition~\ref{cantorconnected} and from Corollary
\ref{cor:basic_divisibility}.
\end{proof}

In particular, we have:

\begin{thm}\label{thm:Urdis}
The Urysohn space $\mathbb{U}_{\Q,\leq 1}$ is divisible.
\end{thm}

In the remainder of this section, we investigate certain conditions
which guarantee that the bounded Urysohn space
$\mathbb{U}_{\Q,\leq 1}$ isometrically embeds into ``large''
parts of itself. These give some measure of the indivisibility of the
space.

We first wish to extend the notion of an orbit and its socket. Indeed
notice that if $S$ is an orbit of the Urysohn space $\mathbb{U}_{V}
=(U;d)$ with socket $B=\{b_i\mid i\in n\}$ and $s$ and $t$ elements in
$S$, then $d(b_i,s)=d(b_i,t)$ for all $i\in n$ because there exists an
isometry which fixes $B$ element-wise and maps $s$ to $t$. Hence we
are led to the following definition.

\begin{defin}\label{defin:socket}
A {\em distance-socket} (or simply {\em d-socket}) $\mathfrak{B}$ of
 $\mathbb{U}_{V}$ is a sequence of the form
\[
\langle (b_0,d_0),(b_1,d_1),\dots,(b_{n-1},d_{n-1})\rangle
\]
so that for all $i,j\in n$:
\begin{enumerate}
\item  $b_i\in U$ and $d_i\in V$. 
\item $d_i+d_j\geq d(b_i,b_j)$.
\item $d_i+d(b_i,b_j)\geq d_j$.
\end{enumerate}
The set $\mathrm{vert}(\mathfrak{B})$ of {\em vertices} of
$\mathfrak{B}$ is the set $\{b_0,b_1,\dots, b_{n-1}\}$, and the set of
{\em distances} of $\mathfrak{B}$ is the set
$\{d_0,d_1,\dots,d_{n-1}\}$.
\end{defin}

An orbit therefore naturaly gives rise to a corresponding socket and
d-socket. But it also follows that given a d-socket
\[
\mathfrak{B}=\langle (b_0,d_0),(b_1,d_1),\dots,(b_{n-1},d_{n-1})\rangle ,
\]
the set of all $s\in U$ so that $d(s,b_i)=d_i$ for all $s\in S$ and
$i\in n$ is an orbit of $\mathbb{U}_V$ with socket $B=\{b_i\mid i\in
n\}$ (Conditions 1. 2. and 3. of the definition ensure that the set is
not empty).

We first show that under certain conditions an orbit itself contains
an isometric copy of the bounded Urysohn space.

\begin{lem}\label{lem:orbchar}
Let $S$ be an orbit of the Urysohn space $\mathbb{U}_{V}=(U;d)$ with
corresponding d-socket $\mathfrak{B}=\langle
(b_0,d_0),(b_1,d_1),\dots,(b_{n-1},d_{n-1})\rangle$.

If $\ell:=\min\{d_i \mid i\in n\}$, then the metric subspace of
$\mathbb{U}_{V}$ induced by $S$ is an isometric copy of the Urysohn
space $\mathbb{U}_{V, \leq 2\ell}$.
\end{lem}

\begin{proof}
Let $i\in n$ be such that $\ell=d_i$. Then $d(s,b_i)=\ell$ for every element
$s\in S$ and hence it follows from the triangle inequality that
$d(s,t)\leq 2\ell$ for any two elements $s$ and $t$ of $S$.

Let $\mathrm{F}:=(F,d^\prime)$ be an element in the age of
$\mathbb{U}_{V, \leq 2\ell}$ so that $F\cap U\subseteq S$ and the
metric subspace of $\mathbb{U}_{V, \leq 2\ell}$ induced by $F\cap
S$ is equal to the metric subspace of $\mathrm{F}$ induced by $F\cap
S$. According to the mapping extension property, it suffices to show that
there exists an embedding of $\mathrm{F}$ into $S$.

Let $\mathrm{G}:=(\{b_i\mid i\in n\}\cup F, \overline{d})$ be the
metric space for which $\overline{d}$ agrees with $d$ on $F \cap U$
and $\overline{d}$ agrees with $d^\prime$ on $F \setminus U$, and
$\overline{d}(x,b_i)=d_i$ for all $x\in F\setminus S$ and all $i\in
n$. The function $\overline{d}$ satisfies the triangle inequality and
hence $\mathrm{G}$ is an element of the age of $\mathbb{U}_{V}$.

\noindent It follows from the mapping extension property of
$\mathbb{U}_{V}$ that there exists an embedding $f$ of $G$
into $U$ which fixes the elements of $G\cap U$. It follows from the
condition that $\overline{d}(x,b_i)=d_i$ for all $x\in F\setminus S$
and all $i\in n$, that the elements of $F$ are mapped by $f$ into $S$.
\end{proof}

If $V'$ is an initial segment of $V$ then $\mathbb U_{V'}$ embeds into
$\mathbb U_V$.  Hence, if we compare orbits of an Urysohn space by
isometric embedding, it follows from Lemma~\ref{lem:orbchar} above
that they form a chain. This is important as you may recall (see
\cite{ES} and \cite{Sauer}) that a necessary condition for
indivisbility is that the ages of the orbits of an homogeneous structure
$\mathrm {H}$ form a chain.  

\begin{cor}\label{norberttest} 
The orbits of an Urysohn metric space form a chain.
\end{cor}

\subsection{Semi-scattered spaces}

We  suppose that $0$ is an accumulation point of $V$.
We show that there are certain small subsets of the Urysohn space $\U_V$  that can be
avoided by any isometrical embedding.

Let $W$ be a subset of $V$ such that $0$ is an accumulation point of $W$.
Thus, in particular, $\mathcal M_W$ contains members of abitrarily  small diameter.

\begin{defin}\label{defin:V-sub-isolation}
An element  $a$ of $\M:=(M,d)\in \mathcal M_V$ is  a {\em
$W$-sub-isolated point } if for every $\epsilon>0$,
there is some non-trivial member $\X\in \mathcal M_W$  of diameter at most $\epsilon$
such that $a$ does not belong  to the union of the isometric copies of $\X$ in $M$. Let $M'$ be the set of elements of $M$ which are not  $W$-sub-isolated and  $\M'$ be the restriction of $\M$ to $M'$. 
\end{defin}

Clearly every isolated point  is $W$-sub-isolated. Thus the following decomposition  generalizes the Cantor-Bendixson decomposition of scattered spaces.

\begin{defin}\label{def:V-sub-scatterness}
Given a  metric space $\M$,  define for each ordinal $\alpha$ 
a metric space $\M^{(\alpha)}$ by
\begin{enumerate}
  \item $\M^{(0)}:=\M$.
  \item $\M^{(\alpha)}:=(\M^{(\beta)})'$ if $\alpha:=\beta+1$.
  \item $\M^{(\alpha)}:=\bigcap_{\beta<\alpha}\M^{(\beta)}$. if $\alpha$ is a limit ordinal.
\end{enumerate}
Clearly, $\beta\leq\alpha\Rightarrow \M^{(\alpha)}\subseteq \M^{(\beta)}$,
hence this ordinal sequence  is eventually constant. It is is eventually empty, we say that $\M$ is {\em $W$-sub-scattered}.
\end{defin}
In the sequel $W:= V$ (and $V$ is a countable subset of $\R_+$ containing $0$, for which $0$ is an accumulation point, and  satisfying the four-values  condition).  

\begin{thm}\label{lem:semi-scatter} 
For every $V$-sub-scattered subspace of $\mathbb{U}_V$
the complementary subspace is isometric to $\mathbb{U}_V$.
\end{thm}

\begin{proof}
First observe that the notion of $V$-sub-isolation is hereditary, {\em i.e.}
a $V$-sub-isolated point of a metric space is $V$-sub-isolated in any subspace it lies in.
It easily follows that every subspace of a $V$-sub-scattered metric space is also $V$-sub-scattered.
 Since $V$ satisfies the four-values  condition,
then, for  every positive real $\ell$,
the metric space $\mathbb{U}_{V,\leq \ell}$ 
has diameter at most $\ell$, 
and in particular it has 
no $V$-sub-isolated point,
since it is homogeneous and it embeds the non singleton $\mathbb{U}_{V,\leq \ell'}$
for any $\ell'\leq\ell$.

Now given a subspace $\M$ of $\U_V$,
if $\M$ is not isometric to $\U_V$, then it follows from 
Lemma~\ref{lem:orbchar} and Point~\ref{prolong:orbit} of Section~\ref{relational}
that the complementary subspace embeds $\U_{V,\leq\ell}$
for some positve $\ell$, and therefore, 
since  $\U_{V,\leq\ell}$ is not $V$-sub scattered,
that complementary subspace is not $V$-sub scattered either. 
\end{proof}

\begin{example}\label{exple:scatter}
For example, given a subspace $\M$ of such a $\U_V$,
the complementary subspace is isometric to $\U_V$
whenever
\begin{itemize}
   \item $\M$ is topologically scattered.
  \item $\M$ is {\em $V$-semi discrete}, 
  {\em i.e.} for every point $a$ of $\M$, $0$ is an accumulation point of $V\setminus Spec(a,\M)$.
\end{itemize}
\end{example}
%


\subsection{ The case of the  Urysohn space $\mathbb{U}_{\Q,\leq 1}$. }

Our final goal is to show that an isometric embedding can avoid a
set containing elements close to a sequence of relatively far elements.

\begin{lem}\label{lem:rim}
Let $\mathfrak{B}:= \langle
(b_0,d_0),(b_1,d_1),\dots,(b_{n-1},d_{n-1})\rangle$ be a d-socket with
orbit $S$ in $\mathbb{U}_{\Q,\leq 1}$.  Let $a\in U$, $r\in \Q\cap [0,1]$, and $x\in S$
such that  $d(a,x)\leq r$.

\noindent Then if $d(a,b_i)\geq r$ for all $i\in n$, there exists an element
$y\in S$ with $d(a,y)=r$.
\end{lem}
\begin{proof}
It suffices to show that $\mathfrak{B}^\prime:= \langle
(b_0,d_0),(b_1,d_1),\dots,(b_{n-1},d_{n-1})(a,r)\rangle$ is again a d-socket. 

\noindent To fulfill the requirements of Definition
\ref{defin:socket}, it remains to verify that for every $i\in n$ the
following inequalities hold:
\begin{enumerate}
\item $r+d_i\geq d(a,b_i)$.
\item $d_i+d(a,b_i)\geq r$.
\item  $r+d(a,b_i)\geq d_i$.
\end{enumerate}

\noindent But by the triangle inequality with $x$ we have
$d(a,b_i)\leq d(a,x)+d_i\leq r+d_i$. Trivially we have $r\leq
d(a,b_i)\leq d_i+d(a,b_i)$, and finally $d_i\leq d(a,x)+d(a,b_i)\leq
r+d(a,b_i)$.
\end{proof}

\begin{thm}\label{thm:holes}
Let $\mathbb{U}_{\Q,\leq 1}=(U;d)$ be the bounded Urysohn
space. Fix $R:=\{r_i \mid i\in \omega\}$ a set of rationals in the
interval $(0,1]$, $A:=\{a_i \mid i\in \omega\}$ a subset of $U$ so
that $d(a_i,a_j)\geq r_i+r_j$ for all $i,j\in \omega$ with $i\not=
j$, and let 
\[
X:=\bigcup_{i\in \omega}\{x\ : \ d(a_i,x)<r_i\}.
\]
Then the metric subspace of $\mathbb{U}_{\Q,\leq 1}$ induced
by $U\setminus X$ is an isometric copy of $\mathbb{U}_{\Q,\leq
1}$.
\end{thm}
\begin{proof}
Notice that if $y$ is any element at distance from some $a_i$ greater
than or equal to $r_i$ for some $i\in \omega$, then $y\in U\setminus
X$.

\noindent Let $S$ be an orbit of $\mathbb{U}_{\Q,\leq 1}$ with
socket $F\subseteq U\setminus X$. 
Given Point~\ref{prolong:orbit} of Section~\ref{relational}, let us check that $S$ meets $U\setminus X$.
Let $s\in S$. If $s\notin X$
then
there is nothing to prove.  Otherwise there exists an $i\in \omega$ so
that $d(a_i,s)<r_i$. But then it follows from Lemma~\ref{lem:rim} that
there is an element $y\in S$ with $d(a_i,y)=r_i$, and hence $y \in
U \setminus X$. This completes the proof.
\end{proof}

\newpage 

\end{document}